\documentclass[11pt]{article}

\usepackage{amsmath}
\usepackage{amsfonts}
\usepackage{amssymb}
\usepackage{epsfig}

\usepackage{amsmath}
\usepackage{amsfonts}
\usepackage{amssymb}
\usepackage{dsfont}
\usepackage{mathrsfs}
\usepackage{bbm}
 \usepackage{hyperref}
\usepackage{amsmath}
\usepackage{amsfonts}
\usepackage{amssymb}
\usepackage{multirow}
\usepackage{mathrsfs}
\usepackage[dvipsnames,usenames]{color}

\renewcommand{\Box}{\framebox{\rule{0.3em}{0.0em}}}

\newtheorem{thm}{Theorem}[section]
\newtheorem{lema}{Lemma}[section]

\newtheorem{ex}{Example}[section]
\newtheorem{rem}{Remark}[section]
\newtheorem{defi}{Definition}[section]

\newtheorem{assu}{Assumption}[section]

\newcommand{\setd}{{ d \kern -.15em l}}
\newcommand{\hatsetd}{ d \hat{\kern -.15em l }}
\newcommand{\dd}{\mathsf {d\kern -0.07em l}} 

\newcommand{\vt}{{\vartheta}}

\newcommand{\bgeqn}{\begin{eqnarray}}
\newcommand{\edeqn}{\end{eqnarray}}
\newcommand{\bgeq}{\begin{eqnarray*}}
\newcommand{\edeq}{\end{eqnarray*}}
\newcommand{\bec}{\begin{center}}
\newcommand{\enc}{\end{center}}
\newcommand{\R}{{\rm I\!R}}

\newcommand{\inmat}[1]{\mbox{\rm {#1}}}

\newcommand{\F}{{\cal F}}

\newcommand{\B}{{\cal B}}

\newcommand{\M}{{\cal M}}

\newcommand{\be}{\begin{equation}}
\newcommand{\ee}{\end{equation}}

\def\e{\epsilon}

\def\bbe{{\Bbb{E}}} 

\renewcommand{\Box}{\hfill \rule{2.3mm}{2.3mm}}

\title{
Statistical Robustness of
Empirical Risks in Machine Learning
}

\author{
Shaoyan Guo\footnote{School of Mathematical Sciences, Dalian University of Technology, Dalian, 116024, China. Email: syguo@dlut.edu.cn.
},
\,\,
Huifu Xu\footnote{Department of Systems Engineering and Engineering Management,The Chinese University of Hong Kong, Shatin, N.T., Hong Kong. Email: hfxu@se.cuhk.edu.hk}
 \,\,and\,\,
Liwei Zhang \footnote{
              School of Mathematical Sciences, 
              Dalian University of Technology, 
              Dalian 116024, China.
              Email: lwzhang@dlut.edu.cn} 
}

\begin{document}

\maketitle

\noindent
\begin{abstract}
This paper studies convergence of empirical risks 
in reproducing kernel Hilbert spaces (RKHS).
A conventional assumption in the existing research is that empirical training data do not contain any noise 
but this may not be satisfied in some practical circumstances.
Consequently 
the  existing convergence results do not provide a guarantee as to 
whether empirical risks
based on empirical data are reliable or not when the data contain some noise.
In this paper, we 
fill out the gap in a few steps.
First, 
we derive moderate sufficient conditions under which the expected risk 
 changes stably (continuously) against small perturbation of the probability distribution of 
 the underlying random  variables and demonstrate how the cost function and kernel
 affect the stability. 
 Second, 
 we examine    the difference between 
    laws of the statistical estimators of the expected 
    optimal loss based on pure data and contaminated 
    data using Prokhorov metric and Kantorovich metric and derive 
    some qualitative and quantitative statistical robustness results.
    Third, we identify appropriate metrics under which the statistical estimators 
    are uniformly asymptotically consistent.
These results provide theoretical grounding for analysing asymptotic convergence and examining reliability of the 
statistical estimators in a number of  well-known machine learning models.
\end{abstract}



\section{Introduction}

%


A key element of 
supervised learning
is to find a function which optimally fits to a training set of input-output data and validate its performance with new data. Classical regression models 
and classification models
are typical examples.
However, with rapid development of social and economic activities and computer technology, data size increases at an exponential rate. This in turn requires much more powerful optimization models 
to understand the behavior of  complex systems with uncertainties on high dimensional parameter spaces
and efficient computational algorithms to solve them. Empirical risk 
minimization
(ERM) is one of them.
The essence of ERM models is to use various approximation methods such as sample average approximation (SAA) and stochastic approximation  
to approximate the expected value of a random function with 
sampled data. 
Regularization is often needed since these problems are usually ill-conditioned. 
Convergence analysis of SAA is well documented in the literature of 
stochastic optimization, see for instance Ruszczy\'nski and Shapiro \cite{RuS03} and references therein.

In the context of machine learning, the focus is not only on the convergence of statistical estimators to their true counterparts as sample size increases, but also on 
scalability of the learning algorithms
because 
the size of machine learning problems 
are often very large
under some circumstances \cite{SSS10}. 
For instance, Norkin and Keyzer \cite{NoK09} consider a general nonparametric regression in finite dimensional 
RKHS and derive nonasymptotic bounds on the minimization
error, exponential bounds on the tail distribution of errors, and sufficient conditions
for uniform convergence of kernel estimators to the true (normal) solution with probability
one. 
In the regularized empirical least squares risk minimization, 
the convergence of estimators can be referred to \cite{CuS02,CuZ07,PoS03,SY06}.
Caponnetto and Vito \cite{CaV07}
develop a theoretical analysis of the performance of the regularized
least-square algorithm 
in the regression setting when the output space is a general Hilbert space.
They use  the concept of effective dimension to choose
the regularization parameter as a function of the number
of samples and derive optimal convergence rates over a suitable class of priors defined by the considered kernel.
More recently,  Davis and Drusvyatskiy \cite{DaD18}
consider a stochastic optimization problem of minimizing population risk,
where the loss defining the risk is assumed to be weakly convex. They establish dimension-dependent rates
on subgradient estimation in full generality and dimension-independent rates when the
loss is a generalized linear model.
We refer readers to  
monograph \cite{CuZ07} for the ML models in infinite dimensional spaces
for a comprehensive overview.

{\color{black}
The problem of characterizing learnability is the most basic question of statistical learning theory. 
For the case of supervised classification and regression,
the learnability is equivalent to uniform convergence of the empirical risk to the expected risk \cite{ABCH97,BEHW89}. For the general 
learning setting, Shalev-Shwartz et al. \cite{SSS10} establish 
that the stability is the key necessary and sufficient condition for 
learnability. The existing literature on stability in learning uses 
many different stability measures. Much of them consider the effect
on the optimal value when there exist small
changes to the sample such as replacing, adding or removing one instant from the sample,
see the review paper \cite{SSS10} for more detail.
A conventional assumption in the above stability is that 
all the instants used in the sample are independent and identically distributed (i.i.d.) and are drawn from the true probability distribution, 
}
but
this may not be satisfied in some practical circumstances,
which means that the empirical training data may contain some noise.
Consequently the  existing convergence results do not provide a guarantee as to whether 
empirical risks and kernel learning estimators
    obtained from solving the ERM models 
    is reliable when the empirical data contain some noise.
    In this paper, we investigate the issue for learning algorithms on a RKHS  
    from statistical robustness perspective \cite{CDS10,KSZ14} in three main steps.

First,  we carry out stability analysis on the optimal expected 
risk of  a generic expected loss minimization problem 
  with respect to perturbation of the probability distribution of the underlying random data. This kind of analysis is well known in stochastic programming (see \cite{Rom03} and references therein)  but not known in 
  machine learning as far as we are concerned. The main challenge in the latter is that the decision variable is often a functional (a function of the underlying random data). In the case when the support of the random data is unbounded, the tail of the probability distribution of the random variables, the tail of the kernel and the tail of the cost function interact and have a joint effect on the stability of the optimal expected risk. We derive moderate sufficient conditions under which the expected risk 
    changes stably (continuously) against small perturbation of the probability distribution
    and demonstrate how the cost function, the kernel and the random data
interactively affect the stability. 

Second, we 
investigate the quality of empirical risk by examining  
the difference between 
    laws of the statistical estimators of the expected 
    risk based on pure data and contaminated 
    data using metrics on probability measures/distributions.
    This kind of approach stems from statistics \cite{Ham71,Hub81,HuR09} and is
    applied to 
  risk management where empirical data are used to estimate risk measures of some random losses by
  Cont et al \cite{CDS10}, Kr\"atschmer et al. \cite{KSZ14,KSZ12} and many others.   
    Here we extend the research to machine learning as we believe the approach  
    can be effectively used to look into the interactions between model errors and data errors 
    from statistical point of view, and we do so 
    in both qualitative and quantitative manners. 

    Third, we discuss convergence of empirical risk which has a vast literature in machine learning.
    Our focus in this paper is on a generic expected loss minimization model in an infinite dimensional 
    RKHS which requires us to take a particular caution on the tails of the kernel and the cost function when they are both unbounded. We also look into the uniform convergence of the statistical estimator with respect to a set of empirical distributions 
    generated near the true one and identify appropriate metrics under which the statistical estimators 
    are uniformly asymptotically consistent.
A combination of all of these results provides some new theoretical grounding for analysing asymptotic convergence and examining reliability of the 
statistical estimators in a number of  well-known machine learning models.

The rest of the paper are organized as follows. 
Section 2 sets up the background of the model and statistical robustness, Section  3 presents stability of 
the expected risk against perturbation of the probability distribution, Section 4 details qualitative and quantitative
analysis of statistical robustness and Section 5 gives uniform consistency analysis, Section 6 points out some future research.

\section{Problem statement}
Let $X$
be the input space and
$Y$
 the output space.
The relation between an input $x\in X$ and
an output $y\in Y$ is described by a probability distribution
$P(x,y)$.
Let $Z$ denote the product space $X\times Y$.
For each input $x\in X$, output $y\in Y$ and $z=(x,y)$, let $c(z,f(x))$ denote
 the  loss caused by  the use of $f$ as a model for the unknown
 process producing $y$ from $x$ and 
 $\bbe_P[c(z,f(x))]$ 
 the 
 statistical average of the losses.
If $P$ is known, then the problem of learning is down to 
find an optimal 
model
such that the average loss is minimized, i.e., 
\bgeqn
\min_{f\in \F} R(f) :=\bbe_P[c(z,f(x))]= \int_Z c(z,f(x))P(dz),
\label{eq:ML-Rf-min}
\edeqn
where $\F$ is some functional class to be specified. 
Let $\vt(P)$ denote the optimal value and 
$\F^*(P)$ 
the set of optimal solutions in (\ref{eq:ML-Rf-min}).
By indicating their dependence on $P$, we will investigate the effect of 
a perturbation of $P$ in forthcoming discussions.
Without loss of generality, we assume throughout the paper that $c(z,f(x))$ takes non-negative value. 
In practice, ${\cal F}$, $Z$ and $c(\cdot,\cdot)$ are known to learners. Here we list a few examples \cite{SSS10}.

\begin{itemize}

    \item \textbf{Regression.} Let $Z=X\times Y$ where $X$ and $Y$ are bounded subsets of $\R^n$ and $\R$ respectively, let
    ${\cal F}$ be a set of functions $f:\R^n\to \R$ and $c(z,f(x))=L(f(x)-y)$,
    where $L(\cdot)$ is a 
    loss function. Specific interesting cases include squared loss function
    $L(t)=\frac{1}{2}t^2$, $\epsilon$-insensitive loss function $L(t)=\max\{0,|t|-\epsilon\}$ with $\epsilon>0$,
    hinge loss function $L(t)=\max\{0,1-t\}$ and log-loss function $L(t)=\log (1+e^{-t})$ in various 
    regression and support vector machine models, see \cite{SKE19}. 
    
     \item \textbf{Binary Classification.} Let $Z=X\times \{0,1\}$ and ${\cal F}$ be a set of functions
    $f:X\to \{0,1\}$, let $c(z,f(x)) = \mathbf{1}_{f(x)\neq y}$. Here $c(\cdot,\cdot)$ is a $0-1$ loss function, 
    measuring whether
    $f$ misclassifies the pair $(x,y)$.
    
    \item \textbf{Density estimation.} Let $Z$ be a subset of $\R^n$ and ${\cal F}$ be a set of bounded probability densities 
on $Z$, let $c(z, f(x))=-\log(f(z))$. Here $c(\cdot,\cdot)$ is simply the negative log-likelihood of an instance $z$ according 
to the hypothesis model density $f$.

\end{itemize}

\subsection{Reproducing kernel Hilbert space}
The
nature of functions $f$ 
in (\ref{eq:ML-Rf-min}) 
needs to be specified.
Let ${\cal H}$ denote a class of functions $f:X\to Y$. ${\cal H}$ is called {\em hypotheses space}
if $f$ 
is restricted to ${\cal H}$. 
This is because the choice of ${\cal H}$ is based on
hypotheses of the structure of these functions.

\begin{defi} Let ${\cal H}(X)$ be a Hilbert space of functions with inner product $\langle \cdot, \cdot\rangle$ and
$k(\cdot,\cdot):X\times X \to \R$ be a kernel, that is, there is a feature map
$\Phi:X\to {\cal H}$ such that $k(x,x) = \langle \Phi(x),\Phi(x)\rangle$.
 ${\cal H}(X)$  is said to be a reproducing kernel Hilbert space (RKHS for short)
 if there is a kernel function $k(\cdot,\cdot):
 X \times X\to \R$ such that:
 (a) $k(\cdot,x)\in {\cal H}(X)$ for all $x\in X$ and (b)
$f(x) =\langle f,k(\cdot,x)\rangle$ for all $f\in {\cal H}(X)$ and $x\in X$.
The corresponding scalar product and norm are denoted by $\langle \cdot,\cdot\rangle$ and $\|\cdot\|_k$ respectively.

A kernel: $k:X\times X\to \R$ is said to be {\em symmetric} if $k(x,t)=k(t,x)$ for each $x,t \in X$,
{\em  positive definite symmetric (PDS) } if for any $x_1,\cdots,x_m \in X$ the matrix
$[k(x_i,x_j)]_{ij}\in \R^{m\times m}$ is symmetric positive semidefinite (SPSD). 
A kernel $k$ is 
called {\em Mercer kernel} if it is continuous, symmetric and positive semidefinite.
\end{defi}



Examples of Mercer kernels abound. Here we list some of them. 

\begin{itemize}

\item \textbf{Polynomial kernel}: 
$k(x_1,x_2)=(\gamma \langle x_1,x_2\rangle +1)^d, \forall x_1,x_2\in \R^N,
$
where
$\gamma>0$ is a constant,
$d \in \mathbb{N}$ and $\mathbb{N}$ denotes the set of positive integers.

\item \textbf{Gussian kernel}:  
$k(x_1,x_2)=e^{-\gamma\|x_1-x_2\|_2^2}, \forall x_1,x_2 \in \R^N,
$
where 
$\gamma >0$ is a constant.

\item 
\textbf{Sigmoid kernel}: 
$
k(x_1,x_2) = \inmat{tanh}\left(a\langle x_1, x_2\rangle+b\right), \forall x_1,x_2 \in \R^N,
$
where $a,b>0$ are constants, $\inmat{tanh}(t)=\frac{e^t-e^{-t}}{e^t+e^{-t}}$ is the hyperbolic tangent function.
 
\end{itemize}

Let $k:X\times X\to \R$ be a positive definite symmetric kernel (Mercer kernel).
 Then there exists a Hilbert space ${\cal H}_k(X)$ and a mapping $\Phi:X \to {\cal H}_k(X)$ such that
$$
k(x,x') =\langle \Phi(x),\Phi(x')\rangle, \forall x,x'\in X.
$$
Moreover ${\cal H}_k(X)$ has the reproducing property, see \cite[Theorem 5.2]{MRT12}.
If we let
\bgeq
\F =\left\{
\sum_{i=1}^n\alpha_ik(x_i,\cdot): n\in \mathbb{N}, \alpha_i\in \R, x_i\in X
\right\}
\label{eq:feasi-F}
\edeq
with the inner product 
$$
\left\langle\sum_{i=1}^n\alpha_ik(x_i,\cdot),\sum_{j=1}^n\beta_jk(x_j,\cdot)\right\rangle = \sum_{i,j=1}^n\alpha_i\beta_jk(x_i,x_j),
$$
then $\F$ can be completed into 
the RKHS, 
see \cite{BBL05}.
 Algorithms working with kernels usually perform minimization of a cost function on a ball of the associated 
 RKHS
 of the form
\bgeqn
\F_\sigma =\left\{
\sum_{j=1}^N a_j k(x_j,\cdot): N\in \mathbb{N}, \sum_{i,j=1}^N a_ia_j k(x_i,x_j) \leq \sigma^2, x_1,\cdots,x_N\in  X
\right\}.
\label{eq:F-SAA}
\edeqn

Throughout the paper, we assume that a positive definite symmetric kernel $k(\cdot,\cdot)$ is given and
${\cal H}_k$ is the RKHS associated with $k$. The functional class $\F$
in (\ref{eq:ML-Rf-min}) and (\ref{eq:ML-saa}) is a subset of 
${\cal H}_k$.

%
%
%
%
%
%
%
%
%
%
%

\subsection{Sample average approximation}

In practice, the true probability distribution $P$ is unknown,
but it is possible
to obtain an independent and identically distributed (i.i.d.) samples $\{z^i=(x^i,y^i)\}_{i=1}^N$ generated by $P$,
which is known as training data.
Given the sample, the goal of machine learning is to find a function $f:X\to Y$ such that $f$
solves
\bgeqn
\min_{f\in \F} 
\bbe_{P_N}[c(z,f(x))] :=\frac{1}{N} \sum_{i=1}^N c(z^i,f(x^i)),
\label{eq:ML-saa}
\edeqn
where
\bgeqn
P_N(\cdot) :=\frac{1}{N} \sum_{i=1}^N \mathbbm{1}_{z^i}(\cdot)
\label{eq:emp-prob-P-N}
\edeqn
denotes the empirical probability measure/distribution and
 $\mathbbm{1}_{z^i}(\cdot)$ denotes the Dirac measure at $z^i$.
Let $\vt(P_N)$ denote the optimal value (empirical risk),
$R_{P_N}(f)$ the objective function,
and $\F^*_{P_N}$  the set of optimal solutions of
the sample average approximation problem (\ref{eq:ML-saa}).
Let $f_N(P_N) \in \F^*_{P_N}$ denote an optimal solution of  (\ref{eq:ML-saa}).
Then  $f_N(P_N)$ is  called the estimator and the framework generating $f_N(P_N)$ is called learning algorithm.
Notice  that from sampling point of view,
we may write $\hat{\vt}_N(z^1,\cdots,z^N)$ and $\hat{f}_N(z^1,\cdots,z^N)$ for $\vt(P_N)$ and $f_N(P_N)$
 respectively to indicate their dependence on the samples.

From computationally perspective,
problem (\ref{eq:ML-saa}) is often ill-conditioned. 
The issue can be  addressed by adopting 
a simple Tikhonov regularization approach:
\bgeqn
\vt(P_N,\lambda_N)=\min_{f\in \F} R_{P_N}^{\lambda_N}(f) :=\bbe_{P_N}[c(z,f(x))]
+ \lambda_N \|f\|_k^2,
\label{eq:ML-saa-r}
\edeqn
where $\lambda_N>0$ is a regularization parameter. In 
general $\lambda_N$ is driven to $0$ but the choice of 
the value may affect the rate of convergence. A  number of papers
have been devoted to this, see for instance 
Brehney and Huang \cite{BrH15} for logistic regression models
in a finite dimensional space, Cucker and Smale 
\cite{CuS02} and Caponnetto and Vito \cite{CaV07} 
for regularized least squares models in infinite dimensional 
RKHS.
Let $\F^*_{P_N,\lambda_N}$ denote the set of optimal solutions in (\ref{eq:ML-saa-r}) and
$f_N(P_N,\lambda_N) \in \F^*_{P_N,\lambda_N}$ 
an
optimal solution. In the case that $c(z,f(x))$ is convex in $f$ for almost every $z$,
$\F^*_{P_N,\lambda_N}$ is a singleton.
By virtue of the representer theorem (see Kimedorf and Wahba \cite{KiW70}, \cite[Theorem 4.2]{ScS02}),
problem (\ref{eq:ML-saa-r}) 
has a solution  which takes the following form
$
f_N^{\lambda_N}(x) =\sum_{j=1}^N \alpha_j k(x_j,x)
$
and by the reproducing property (\cite{NoK09}),
$
\|f_N^{\lambda_N}\|_k^2 
= \langle f_N^{\lambda_N},f_N^{\lambda_N}\rangle =\sum_{i,j=1}^N \alpha_i\alpha_j k(x_i,x_j).
$
In that case  the feasible set may be written as (\ref{eq:F-SAA}).
As we commented earlier, here we may also write $\hat{\vt}_N(z^1,\cdots,z^N,\lambda_N)$ and $\hat{f}_N(z^1,\cdots,z^N,\lambda_N)$ for $\vt(P_N,\lambda_N)$ and $f_N(P_N,\lambda_N)$
 respectively to indicate their dependence on the samples.

\subsection{Contamination of the training data}


The current research of machine learning is mostly focused on the case that sample data are 
generated by the true probability distribution $P$ which means that they do not contain any noise.
This assumption may not be 
satisfied in some data-driven problems.
Let $\tilde{z}^1,\cdots,\tilde{z}^N$ denote the 
perceived data which may contain noise
and 
\bgeqn
Q_N(\cdot) :=\frac{1}{N} \sum_{i=1}^N \mathbbm{1}_{\tilde{z}^i}(\cdot)
\label{eq:emp-prob-Q-N}
\edeqn
be the respective empirical distribution.
Instead of solving (\ref{eq:ML-saa-r}), 
we 
solve, in practice,
\bgeqn
\min_{f\in \F}
\bbe_{Q_N}[c(z,f(x))]
+ \lambda_N \|f\|_k^2.
\label{eq:ML-saa-r-Q}
\edeqn
Let $R_{Q_N}^{\lambda_N}(f)$, $\vt(Q_N,\lambda_N)$ and $f_N(Q_N,\lambda_N)$
denote respectively the objective function, the 
optimal value and the optimal solution of problem (\ref{eq:ML-saa-r-Q}).
We 
are then concerned with the quality of the learning model estimator $f_N(Q_N,\lambda_N)$ and the associated empirical risk  $\vt(Q_N,\lambda_N)$. 

There are two ways to proceed the research.
One is to look into convergence of the statistical quantities as the sample size $N$ increases
and the regularization parameter $\lambda_N$ 
goes to  zero. 
Assume without loss of generality that the samples are independent and identically generated.
By law of large numbers, $Q_N$ converges to some probability distribution $Q$ and subsequently  
\bgeqn
f_N(Q_N,\lambda_N) \to f(Q)  \quad \inmat{and} \quad  \vt(Q_N,\lambda_N) \to \vt(Q).
\label{eq:asyconvg}
\edeqn
On the other hand, if we regard $Q$ as a perturbation of the true unknown probability distribution $P$, 
then we need to investigate whether
\bgeqn
f(Q)\to f(P)  \quad \inmat{and} \quad \vt(Q)\to\vt(P)
\label{eq:continuity}
\edeqn
as $Q$ approaches $P$. 
The former is known as asymptotic convergence/consistency 
 and the latter is known as stability in the literature of stochastic programming \cite{Rom03}. 
However, if we want to establish 
\bgeqn
f_N(Q_N,\lambda_N) \to f(P)  \quad \inmat{and} \quad  \vt(Q_N,\lambda_N) \to \vt(P),
\label{eq:stability} 
\edeqn
then we require not only (\ref{eq:continuity}) but also
(\ref{eq:asyconvg}) to hold uniformly for all $Q$ near $P$.
This will be more demanding than the currently established convergence results.

The other is to examine the 
discrepancy 
between $f_N(Q_N,\lambda_N)$
and $f_N(P_N,\lambda_N)$ ($\vt(Q_N,\lambda_N)$
and $\vt(P_N,\lambda_N)$) via law of these estimators. 
The latter should be understood as estimators when the noise in the samples is detached (an ideal case).  
This kind of research 
is in alignment with qualitative robustness in the literature of robust statistics and risk measurement,
see \cite{CDS10,GuX18,KSZ14,KSZ12} and references therein.  We will give a formal definition
in Section 4.

In both steps leading towards statistical robustness of $\vt(\cdot)$, we will need to restrict
the perturbation of the probability measure from $P$ 
to 
the space of $\phi$-topology of weak convergence 
instead of usual
weak convergence.


\subsection{$\phi$-weak topology}

We recall some basic concepts and results
about weak topology which are needed for the analysis. The materials are
  mainly  extracted  from \cite{Clau16}, we refer readers to 
  \cite[Chapter 2]{Clau16} and references therein
  for a more comprehensive discussion on the subject.

\begin{defi} Let $\phi:Z\to [0,\infty)$ be a continuous function and
\bgeq
{\cal M}_Z^\phi := \left\{P\in \mathscr{P}(Z): \int_{Z} \phi(z)P(dz)<\infty\right\},
\edeq
where
$\mathscr{P}(Z)$ is the set of all probability measures on the measurable space $\left(Z,\B(Z)\right)$ with Borel sigma algebra $\mathcal B(Z)$ of $Z$.
\end{defi}

${\cal M}_Z^\phi $ defines a subset of probability measures in $\mathscr{P}(Z)$ which satisfies the generalized moment condition
of $\phi$.

\begin{defi}[$\phi$-weak topology] Let $\phi:Z\to [0,\infty)$ be a gauge function, that is, $\phi\geq 1$ holds outside a compact set.
Define ${\cal C}_Z^\phi$ the linear space of all continuous functions $h:Z\to \R$ for which there exists a positive constant $c$ such that
$$
h(z)\leq c(\phi(z)+1), \forall z\in Z.
$$
The $\phi$-weak topology, denoted by $\tau_\phi$, is the coarsest topology on
${\cal M}_Z^\phi$ for which the mapping $g_h:{\cal M}_Z^\phi\to \R$ defined by
$$
g_h(P) :=\int_{Z} h(z) P(dz),\; h\in {\cal C}_Z^\phi
$$
 is continuous. A sequence  $\{P_l\} \subset {\cal M}_Z^\phi$ is said to converge $\phi$-weakly to $P\in {\cal M}_Z^\phi$ written
 ${P_l} \xrightarrow[]{\phi} P$ if it converges 
 with respect to (w.r.t.) $\tau_\phi$.
\end{defi}

From the definition, we can see immediately that $\phi$-weak convergence implies weak convergence under usual topology of weak convergence.
We denote the latter by ${P_l} \xrightarrow[]{w} P$.
 Moreover, it follows by \cite[Corollary 2.62]{Clau16} that the $\phi$-weak topology on ${\cal M}_Z^\phi$ is generated by the metric
$\dd_\phi:{\cal M}_Z^\phi\times {\cal M}_Z^\phi\to \R$ defined by
\bgeqn
\dd_\phi(P',P''):=\dd_{\inmat{Prok}}(P',P'')+\left|\int_{Z} \phi d P'-\int_{Z} \phi d P''\right|, \; \inmat{for} \; P', P'' \in {\cal M}_Z^{\phi},
\label{eq:d-psi}
\edeqn
where $\dd_{\inmat{Prok}}: \mathscr{P}(Z)\times \mathscr{P}(Z)\to \R_+$ is the Prokhorov metric defined as follows:
\bgeqn
\dd_{\inmat{Prok}}(P',P''):=\inf\{
\epsilon>0: P'(A) \leq P''(A^\epsilon)+\epsilon\,\mbox{for all}\,A \in \mathcal B(Z)
\},
\edeqn
where $A^\epsilon:= A + B_\epsilon(0)$ denotes the Minkowski sum of $A$ and the open ball centred at $0$ (w.r.t. the norm in $Z$).
When $\phi \equiv 1$,  the second term in (\ref{eq:d-psi}) disappears and consequently $d_\phi(P',P'')=d_{\inmat{Prok}}(P',P'')$.
In that case, the $\phi$-weak topology reduces to the usual topology of weak convergence (defined through bounded continuous functions).
Equivalence between the two topologies may be established over a set which satisfies some uniform integration conditions,
see \cite[Lemma 2.66]{Clau16} and the reference therein.

\begin{defi}[Fortet-Mourier metric]
\label{D-Fort-Mou-metric}
Let 
$$
\mathcal{F}_{p}(Z):=\left\{\psi: Z\rightarrow \R: |\psi(z)-\psi(\tilde{z})|\leq c_{p}(z,\tilde{z})\|z-\tilde{z}\|, \forall z,\tilde{z}\in Z\right\},
$$
where $\|\cdot\|$ denotes some norm on $Z$ and $c_{p}(z,\tilde{z}):=\max\{1,\|z\|,\|\tilde{z}\|\}^{p-1}$ 
for all $z,\tilde{z}\in Z$ and $p\geq 1$ describes the growth of the local Lipschitz constants. 
The $p$-th order Fortet-Mourier metric over $\mathscr{P}(Z)$ is defined by
\bgeqn
\zeta_{p}(P,Q):=\sup_{\psi\in \mathcal{F}_{p}(Z)}\left|\int_{Z}\psi(z)P(dz)-\int_{Z}\psi(z)Q(dz)\right|.
\label{eq:zetametric}
\edeqn
\end{defi}

Fortet-Mourier metric is well-known in stochastic programming. 
The unique feature of the metric is that it is induced by a class of locally Lipschitz continuous functions with
specified modulus and rate of growth.  In the case when $p=1$, it reduces to 
Kontorovich metric.  We refer readers to see R\"omisch \cite{Rom03} for a comprehensive overview of the topic.  
From the definition, we can see that 
$$
\zeta_p(P,Q)\leq \bbe_{P\times Q}[c_p(z,\tilde{z})\|z-\tilde{z}\|],
$$ 
where
$P\times Q$ denotes the joint probability distribution of $z$ and $\tilde{z}$. In the case when $P$ and $Q$
are empirical distributions generated by i.i.d. sample, we have
$$
\bbe_{P\times Q}[c_p(z,\tilde{z})\|z-\tilde{z}\|] = \frac{1}{N^2}\sum_{i,j=1}^N c_p(z^i,\tilde{z}^j)\|z^i-\tilde{z}^j\|.
$$
The latter may be used to give an estimate of $\zeta_p(P,Q)$ if we are able to obtain the i.i.d. samples in practice.

\section{Stability analysis} 

In this section, 
we investigate how the model risk of problem (\ref{eq:ML-Rf-min}) is affected by a small perturbation of the probability measure $P$. 
This kind of research is well known in the literature of stochastic programming \cite{Rom03} but
not in machine learning as far as we are concerned. 
We proceed with some technical assumptions which stipulate the properties of the cost function and the kernel.

\begin{assu}
\label{A:kernel}
\begin{itemize}
\item[(a)]For any compact subset $Z_0$ of $Z$, let $X_0$ be its orthogonal  projection 
on $X$. The set of functions $\{k(\cdot,x): x\in X_0\}$ are
equi-continuous on $X_0$, i.e.,
for any $\epsilon>0$, there exists a constant $\eta>0$ such that
$$
\|k(\cdot,x')-k(\cdot,x)\|_k < \epsilon, \forall x, x' \in X_0:  \|x'-x\| < \eta,
$$
 where $\|\cdot\|$ is 
 some norm on $X$.

\item[(b)] There is a positive constant $\beta$ such that  $\|f\|_k \leq \beta$ for all $f\in \F$.
 \end{itemize}
\end{assu}


\begin{rem}\label{r:calm}
To see how 
Assumption \ref{A:kernel} (a) can be possibly satisfied, we recall the notion of 
calmness of kernel introduced by Shafieezadeh-Abadeh et al. \cite[Assumption 25]{SKE19}.
The kernel function $k$ is said to be  {\em calm from above},
if there exists a concave smooth growth function $g:\R_+ \to \R_+$ with $g(0)=0$ and 
$g'(t) \geq 1$ for all $t\in \R_+$ such that
\bgeq
\sqrt{k(x',x')-2k(x,x')+k(x,x)}
\leq g(\|x-x'\|), \forall x, x \in X.
\edeq
Under the calmness condition, there exists $\eta>0$ such that
\bgeq
\|k(\cdot,x')-k(\cdot,x)\|_k
&=&\sqrt{\langle k(\cdot,x')-k(\cdot,x),k(\cdot,x')-k(\cdot,x) \rangle}\\
&=&\sqrt{ k(x',x')-k(x,x')+k(x,x)-k(x,x') }\\
&\leq& g(\|x-x'\|) < \epsilon
\edeq
for all $x,x'$ with $\|x-x'\|\leq \eta$.
The last inequality is due to the fact that 
the growth function $g$ is continuous with
$g(0)=0$, thus for any $\epsilon >0$, there exists a positive constant $\eta >0$ such that $|g(t)-g(0)|=|g(t)|<\epsilon$
for $|t| < \eta$. 
The calmness condition is non-restrictive, which can be satisfied in the following 
cases for $X=\R^n$, see \cite[Example 1]{SKE19}.

\begin{itemize}
\item Linear kernel: for $k(x_1,x_2)=\langle x_1,x_2\rangle$,
$g(t)=t$.

\item Gaussian kernel: for $k(x_1,x_2)=e^{-\gamma\|x_1-x_2\|_2^2}$, 
$g(t)=\max\{\sqrt{2\gamma},1\}t$.

\item Laplacian kernel: for $k(x_1,x_2)=
e^{-\gamma \|x_1-x_2\|_1}$, 
$g(t)= \sqrt{2\gamma t \sqrt{n}}$
if $0\leq t \leq \gamma \sqrt{n}\slash 2$
and $g(t)=t+\gamma \sqrt{n}\slash 2$ otherwise.

\item Polynominal kernel: the kernel $k(x_1,x_2)=(\gamma \langle x_1,x_2\rangle +1)^d$
with $\gamma>0$ and $d \in \mathbb{N}$
fails to satisfy the calmness condition if
$X$ is unbounded and $d>1$, in which case
$\sqrt{k(x_1,x_1)-2k(x_1,x_2)+k(x_2,x_2)}$
grows superlinearly.
If $X \subset \{x\in \R^n: \|x\|_2 \leq R\}$
for some $R >0$, however,
the polynomial kernel is calm with respect to the
growth function
\bgeq
g(t)=
\left\{
\begin{array}{ll}
\max\{\frac{1}{2R}\sqrt{2(\gamma R^2 +1)^d},1\}t &d\, \inmat{is even},\\
\max\{\frac{1}{2R}\sqrt{2(\gamma R^2 +1)^d}-
2(1-\gamma R^2)^d,1\}t &d\, \inmat{is odd}.
\end{array}
\right.
\edeq
\end{itemize}

Assumption \ref{A:kernel} (b) may be guaranteed by restricting the set of feasible solutions to lie within a ball, see
\cite[Assumption D]{NoK09}.
\end{rem}

\begin{assu}
\label{A:cost-1}
The cost function $c(\cdot,\cdot)$ satisfies the following properties.
\begin{itemize}

\item[(a)]
There is a
gauge function $\phi(\cdot)$ such that
\bgeqn
c(z,f(x))\leq \phi(z), \forall z\in Z \; \inmat{and} \; f\in \F,
\label{eq:growth-ml}
\edeqn
where $\phi(z)\to \infty$ as $\|z\|\to \infty$.

\item[(b)] $c(z,y): Z\times Y\to \R$ is continuous.
\end{itemize}

\end{assu}

\begin{rem} 
Condition (a) is known as a growth condition where $\phi(z)$ controls the growth of the cost function as $\|z\|$ goes to infinity.
It 
is trivially satisfied when $Z$ is compact. Our focus here is on the case that $Z$ is unbounded.
  Obviously $\phi$ depends on the concrete structure of $c(.,.)$. 
Consider for example $c(z,f(x))=\frac{1}{2}\|y-f(x)\|^2$. Then
\bgeq
c(z,f(x))
&\leq& \|y\|^2+ \|f(x)\|^2
=  \|y\|^2+ |\langle f, k(\cdot,x)\rangle|^2\\
&\leq& \|y\|^2+ \|f\|_k^2 \|k(\cdot,x)\|_k^2.
\edeq
Moreover, under Assumption \ref{A:kernel} (b), i.e., $\|f\|_k \leq \beta$,
we can work out an explicit form of $\phi$ for some specific kernels.
\begin{itemize}
\item If $k$ is a Linear kernel, then
$\|k(\cdot,x)\|_k^2=|k(x,x)|=\|x\|^2$ and 
$ 
\phi(z):=\|y\|^2+\beta^2\|x\|^2;
$

\item If $k$ is a Gaussian kernel or Laplacian kernel, then 
$\|k(\cdot,x)\|_k^2=0$ and 
$
\phi(z):=\|y\|^2.$

\item If $k$ is a Polynominal kernel, then 
$\|k(\cdot,x)\|_k^2=(\gamma \|x\|^2+1)^{d}$ and 
\bgeqn
\phi(z):=\|y\|^2+\beta^2(\gamma \|x\|^2+1)^{d}.
\label{eq:Poly-kel-phi}
\edeqn
\end{itemize}
From the examples above, we can see that $\phi$ captures not only the growth of the cost function $c(\cdot,\cdot)$ but also the kernel. The growth rate of $\phi$ at the tail in turn affects the topology of weak convergence to be used in the stability analysis in the next theorem.

\end{rem}

\begin{thm}
\label{t-cont-vt-wrt-P}
Under Assumptions \ref{A:kernel} and \ref{A:cost-1}, the following holds for any $p \geq 1$,
\bgeqn
  \label{eq:conti-vt}
\lim_{
P'\xrightarrow[]{\phi^p} P} \vt(P')=\vt(P).
\edeqn
\end{thm}

\noindent
\textbf{Proof.}
Since $(\M^{\phi^p}_Z, \tau_{\phi^p})$ is a Polish space,
it suffices to show that (\ref{eq:conti-vt}) holds for any sequence
 $\{P_l\}\subset \M^{\phi^p}_Z$
with $P_l \xrightarrow[]{\phi^p} P\in \M^{\phi^p}_Z$.
First, 
$P_l \xrightarrow[]{\phi^p} P$ implies that
$P_l \xrightarrow[]{w} P$  and
$$
\lim_{l \to \infty} \int_Z\phi^p(z) P_l(dz)=\int_Z \phi(z)^p P(dz).
$$
Moreover, by  \cite[Lemma 2.61]{Clau16},  for any $\epsilon>0$, there exists a positive constant $M>1$ such that
\bgeqn
\int_Z\phi^p(z) \mathbbm{1}_{(M,\infty)}(\phi^p(z)) P(dz)< \epsilon
\label{eq:cont-prf-M-a}
\edeqn
and
\bgeqn
\sup_{l\in \mathbb{N}} \int_Z\phi^p(z) \mathbbm{1}_{(M,\infty)}(\phi^p(z)) P_l(dz) < \epsilon,
\label{eq:cont-prf-M-b}
\edeqn
where 
$\mathbbm{1}_{(M,\infty)}(t)=1$ if $t\in (M, \infty)$ otherwise $0$. Since $\phi$ 
is coercive, i.e., $\phi^p(z) \to \infty$ as $\|z\| \to \infty$,
then  exists a compact continuity set $Z_M\subset Z$ of $P$ such that 
$Z \backslash Z_M \subset \{z \in Z: \phi^p(z)>M\}$.
Here the continuity set means that $P(\partial Z_M)=0$ where
$\partial Z_M$ denotes the boundary of $Z_M$. 

Let
$$
\mathscr{G} := \{g:
g(z):=c(z,f(x))\;\inmat{for}\; f \in \F\}
$$
and
$$
\mathscr{G}_{M}:=\{g_{M}: Z_M \to \R|
g_{M}(z):=g(z)\; \inmat{for}\; z \in Z_M, g \in \mathscr{G}\}.
$$
It follows from Assumption \ref{A:cost-1} (a) that for each $g_M \in \mathscr{G}_{M}$
and $z \in Z_M$,
$|g_M(z)|\leq \sup_{z\in Z_M} \phi(z) <\infty$, which implies that $\mathscr{G}_{M}$
is uniformly bounded. 

Next, we prove that $\mathscr{G}_{M}$ is equi-continuous over $Z_M$.
By the reproducing property of the kernel
$k(\cdot,\cdot)$, i.e.,
$f(x) = \langle f, k(\cdot,x)\rangle$ for every $f\in \F$, 
we have
\bgeqn
|f(x')-f(x)|
&=& |\langle f, k(\cdot,x')\rangle-\langle f, k(\cdot,x)\rangle|\leq \|f\|_k\|k(\cdot,x')-k(\cdot,x)\|_k\nonumber\\
&\leq& \beta \|k(\cdot,x')-k(\cdot,x)\|_k.
\label{eq:g-equi-cont-ML-0}
\edeqn
The equicontinuity of $k(\cdot,x)$  over $X_M$ (under  Assumption \ref{A:kernel} (a))
ensures the equicontinuity of $\F$ over $X_M$.
Moreover, 
under Assumption \ref{A:cost-1}(b), $\mathscr{G}_{M}$ is also equicontinuous because
$c(\cdot,\cdot)$ is uniformly continuous over any compact set.

Let $Q_l, Q$ be measures on $Z_M$ defined by 
$Q_l(A)=P_l(A)$ and $Q(A)=P (A)$ respectively. 
Since $Z_M$ is a continuity set of $P$, 
then $P_l \xrightarrow[]{w} P$ imply 
$Q_l \xrightarrow[]{w} Q$.
Since $\mathscr{G}_{M}$ is uniformly bounded and equi-continuous,
by \cite[Theorem 3.1]{Rao62},
 \bgeqn
 \lim_{l\to \infty} \sup_{g_M\in \mathscr{G}_M}
\left| \int_{Z_M} g_{M}(z) Q_l(dz)-
 \int_{Z_M} g_{M}(z) Q(dz) \right|=0.
\label{eq:conti-vt-equiv}
\edeqn
On the other hand, under the growth condition (\ref{eq:growth-ml}), (\ref{eq:cont-prf-M-a})
and (\ref{eq:cont-prf-M-b})
imply
\bgeqn
\int_{Z\backslash Z_M} |g(z)| P(dz)
\leq \int_{Z\backslash Z_M} \phi^p(z)  P(dz)\leq 
\int_Z\phi^p(z) \mathbbm{1}_{(M,\infty)}(\phi^p(z)) P(dz)<\epsilon
 \label{eq:unfm-int-phi-ML-1}
 \edeqn
 and
 \bgeqn
\sup_{l\in \mathbb{N}} \int_{Z\backslash Z_M} |g(z)|  P_l(dz)
\leq \sup_{l\in \mathbb{N}} \int_{Z\backslash Z_M} \phi^p(z)  P_l(dz)
\leq \sup_{l\in \mathbb{N}} \int_Z\phi^p(z) \mathbbm{1}_{(M,\infty)}(\phi^p(z)) P_l(dz)
<\epsilon.\nonumber\\
 \label{eq:unfm-int-phi-ML-2}
\edeqn
Together with (\ref{eq:conti-vt-equiv}), we have
\bgeq
&&|\vt(P_l)-\vt(P)|\\
&\leq& \sup_{f\in \F}
\left| \int_Z c(z,f(x))P_l(dz)-
 \int_Z c(z,f(x))P(dz) \right|\\
 &=&\sup_{g\in \mathscr{G}}
\left| \int_Z g(z)P_l(dz)-
 \int_Z g(z)P(dz) \right|\\
 &\leq& \sup_{g\in \mathscr{G}}
\left| \int_{Z_M} g(z) P_l(dz)-
 \int_{Z_M} g(z) P(dz) \right|+\int_{Z\backslash Z_M} |g(z)| P(dz)
 +\int_{Z\backslash Z_M} |g(z)|  P_l(dz)
 \\
 &\leq&\sup_{g_M\in \mathscr{G}_M}
\left| \int_{Z_M} g_{M}(z) Q_l(dz)-
 \int_{Z_M} g_{M}(z) Q(dz) \right|+2\epsilon <3\epsilon
\edeq
for sufficiently large $l$.
The proof is complete.
\hfill $\Box$

The theorem tells us that $\vt(Q)$ is close to $\vt(P)$
when $Q$ is perturbed from $P$ under the $\tau_{\phi^p}$-weak topolgy for any fixed $p \geq 1$. 
Since the empirical probability measure $P_N\in {\cal M}_Z^{\phi^p}$,
we have 
\bgeqn
  \label{eq:conti-vt-SAA}
\lim_{N\to \infty}  \vt(P_N)=\vt(P)
\edeqn
almost surely.  The topological structure of 
set ${\cal M}_Z^{\phi^p}$ affects the stability of $\vt(\cdot)$:
a larger ${\cal M}_Z^{\phi^p}$ 
means that $\vt(\cdot)$ remains stable w.r.t.  
a greater freedom of perturbation from $P$.
In the case when $Z$ is a compact set, ${\cal M}_Z^{\phi^p}=\mathscr{P}(Z)$, which means 
$\vt(\cdot)$ remains stable for any perturbation of the probability measure from $P$ locally.
The tail behaviour of $c(z,f(x))$ affects 
the structure of  ${\cal M}_Z^{\phi^p}=\mathscr{P}(Z)$, we explain this through next example.

\begin{ex}
Consider the least squares regression model with
Polynomial kernel. By (\ref{eq:Poly-kel-phi})
\bgeq
{\cal M}_Z^\phi &=& \left\{P\in \mathscr{P}(Z): \int_Z \left[\|y\|^2+\beta^2(\gamma \|x\|^2+1)^{d}\right]P(dz)<\infty\right\}\\
&=& \left\{P\in \mathscr{P}(Z): \int_Z \|y\|^2 P(dz)<\infty, 
\int_Z \|x\|^{2d} P(dz)<\infty\right\}.
\edeq
We can see from the formulation above that a larger $d$ requires a thinner tail of $P$ and hence a smaller set of  ${\cal M}_Z^\phi$, consequently the stability result is valid for a smaller class of probability distributions.

In the case of Gaussian kernel or Laplacian kernel,
\bgeq
{\cal M}_Z^\phi = \left\{P\in \mathscr{P}(Z): \int_Z \|y\|^2 P(dz)<\infty\right\},
\edeq
which is the set of probability  measures with finite second order moment of $y$.

\end{ex}

Finally, we note that our stability result should be distinguished from those in \cite{SSS10} where stability 
is used to examine the difference of the costs resulting from kernel learning estimators based on different samples.

\section{Statistical robustness}

We now move on to discuss statistical robustness of the machine learning model (\ref{eq:ML-saa-r}).
To ease the exposition, let $Z^{\otimes N}$ denote  the Cartesian product  $Z \otimes\cdots\otimes Z$ and
$\B(Z)^{\otimes N}$ its Borel sigma algebra. 
Let  $P^{\otimes N}$  denote
the probability measure on the measurable space
$\left(Z^{\otimes N},\B(Z)^{\otimes N}\right)$ with marginal $P$ 
and
$Q^{\otimes N}$ with marginal $Q$. 
We will consider statistical estimators mapping from $\left(Z^{\otimes N},\B(Z)^{\otimes N}\right)$ to $\R$
and examine their convergence under $Q^{\otimes N}$ and
$P^{\otimes N}$.

\subsection{Qualitative robustness}

We begin by
a formal definition of  statistical estimator $T(\cdot,\lambda)$ parameterized by $\lambda$,
where $T(\cdot,\lambda)$
maps from a subset of $\mathcal{M} \subset \mathscr{P}(Z)$ to $\R$.
To ease the exposition, 
we 
write
 $\vec{z}^N$ for $(z^1,\cdots,z^N)$ and 
 $\hat{T}_N(\vec{z}^N,\lambda_N)$
for $T(P_N,\lambda_N)$ for fixed sample size $N$.
The following definition is based on
 Kr\"atschmer et al.  \cite[Definition 2.11]{KSZ14}.

\begin{defi}[Statistical robustness] Let
${\cal M} \subset \mathscr{P}(Z)$ be a set of probability measures
and $\dd_\phi$ be defined as in (\ref{eq:d-psi}) for some gauge function $\phi:Z\to \R$, let $\{\lambda_N\}$ be a sequence of parameters.
 A parameterized  statistical estimator $T(\cdot,\lambda_N)$ is said to be robust on
 $\cal{M}$ with respect to $\dd_\phi$ and $\dd_{\inmat{Prok}}$ if for all
$P \in {\cal M}$ and $\epsilon>0$, there exist $\delta>0$ and $N_0\in \mathbb{N}$ such that for all $N\geq N_0$
\bgeq
Q\in {\cal M}, \dd_\phi(P,Q) \leq \delta \Longrightarrow 
\dd_{\inmat{Prok}}\left(P^{\otimes N}\circ \hat{T}_N(\cdot,\lambda_N)^{-1},
Q^{\otimes N}\circ \hat{T}_N(\cdot,\lambda_N)^{-1}
\right) \leq \epsilon.
\edeq
\end{defi}

In this definition, $P^{\otimes N}\circ \hat{T}_N(\cdot,\lambda_N)^{-1}$
and 
$
Q^{\otimes N}\circ \hat{T}_N(\cdot,\lambda_N)^{-1}
$
are 
two probability distributions of random variable $\hat{T}_N(\cdot,\lambda_N)$
mapping  from 
probability spaces 
$\left(
Z^{\otimes N},\B(Z)^{\otimes N}, P^{\otimes N}
\right)$
and
$\left(
Z^{\otimes N},\B(Z)^{\otimes N}, Q^{\otimes N}
\right)$ respectively
to $\R$, and the Prokohorov metric is used to measure
the difference of the two distributions (also known as laws in the literature \cite{CDS10,KSZ14}). The statistical robustness
requires the difference under the Prokhorov metric
to be small when the difference between $P$ and $Q$ is small under
$\dd_\phi$. The definition relies heavily on the adoption of
the two metrics. In Cont et al. \cite{CDS10}, the authors use 
L\'evy metric for both.
Kr\"atschmer et al. \cite{KSZ14} argue that the L\/evy metric underestimates the impact of the tail distributions of $P$ and $Q$ and subsequently propose to use $\dd_\phi$ to replace the L\'evy metric.
Since the former is tighter than the later, it means the perturbation under $\dd_\phi$ is more restrictive and hence 
 enables one to 
examine finer difference between the laws of the statistical estimators.

Statistical robustness is also called qualitative robustness in this paper in that there is no explicit quantitative  relationship between $\epsilon$ and $\delta$.
To establish the statistical robustness,
we need the following Uniform Glivenko-Cantelli property.

\begin{defi}[Uniform Glivenko-Cantelli property]
{\color{black} Let $\phi$ be a gauge function and $\dd_\phi$ be defined as in (\ref{eq:d-psi}).
Let ${\cal M}$ be a subset of $\mathcal{M}_Z^\phi$.
the metric space
$(\cal M,\dd_\phi)$} 
is said to have 
Uniform Glivenko-Cantelli (UGC) property
if for every $\epsilon>0$ and $\delta>0$, there exists $N_0\in \mathbb{N}$
such that for all $P\in {\cal M}$
\bgeqn
P^{\otimes N}\left[\vec{z}^N:
\dd_\phi(P,P_N) \geq \delta\right] \leq \epsilon, \forall P\in \M
\label{eq:UGC}
\edeqn
for all $N\geq N_0$.
\label{d-UGC}
\end{defi}

Recall that $P_N$ is constructed through i.i.d. samples 
generated by random variable $z$ with probability distribution $P$.
The UGC property requires that 
for all $P\in \M$,
their empirical probability measures converge to 
their true counterparts uniformly as the sample size goes to infinity.
The convergence is under $\dd_\phi$ which means not only the weak convergence but also convergence 
of the $\phi$ moments, the latter captures the tails of $P$.

\begin{thm}[Statistical robustness]
\label{t-SR-u-indentity}
Let $\{P_N\}$ be a sequence of empirical probability measures defined by (\ref{eq:emp-prob-P-N})
 and
$ \M_{Z,\kappa}^{\phi^p}$ be the class of all $P\in \mathscr{P}(Z)$ such that
\bgeqn
\int_{Z} \phi(z)^{p}  P(dz) \leq\kappa,
\label{eq:M-k-kappa}
\edeqn
for $\kappa \geq 0$ and $p>1$.
Let Assumptions~\ref{A:kernel} and \ref{A:cost-1}
hold, 
$\lambda_N\to 0$ as $N\to \infty$.
Then for any $\epsilon>0$, there exist positive numbers
 $\delta>0$ and $N_0\in \mathbb{N}$ such that
when $ Q \in {\cal M} \subset  {\cal M}_{Z,\kappa}^{\phi^p}$, $\dd_{\phi}(P,Q) \leq \delta$,
we have
\bgeqn
\dd_{\inmat{Prok}}\left(P^{\otimes N}\circ \hat{\vt}_N(\cdot,\lambda_N)^{-1},
Q^{\otimes N}\circ \hat{\vt}_N(\cdot,\lambda_N)^{-1}\right) \leq \epsilon
\label{eq:stat-robust-vt_n-SP}
\edeqn
 for all $N\geq N_0$ 
 and $\lambda_N \leq \frac{\epsilon}{6\beta^2}$, where $\hat{\vt}_N(\vec{z}^N,\lambda_N)= \vt(P_N,\lambda_N)$ denotes the optimal value of problem (\ref{eq:ML-saa-r}).

\end{thm}

\noindent\textbf{Proof.}
By triangle inequality
\bgeq
\label{eq:stat-robust-vt_n-SP-1}
 &&\dd_{\inmat{Prok}}\left(P^{\otimes N}\circ \hat{\vt}_N(\cdot,\lambda_N)^{-1},
Q^{\otimes N}\circ \hat{\vt}_N(\cdot,\lambda_N)^{-1}\right) \\
&\leq& 
\dd_{\inmat{Prok}}\left(P^{\otimes N}\circ \hat{\vt}_N(\cdot,\lambda_N)^{-1},
\mathbbm{1}_{\inf_{f\in \F} R_P(f)}\right)+
\dd_{\inmat{Prok}}\left(
\mathbbm{1}_{\inf_{f\in \F} R_P(f)},\mathbbm{1}_{\inf_{f\in \F} R_Q(f)}\right)\nonumber\\
&&+
\dd_{\inmat{Prok}}\left(
\mathbbm{1}_{\inf_{f\in \F} R_Q(f)},
Q^{\otimes N}\circ \hat{\vt}_N(\cdot,\lambda_N)^{-1}
\right),\nonumber
\edeq
where $\mathbbm{1}_{a}$ denotes the Dirac  measure at $a\in \R$. 
By Theorem \ref{t-cont-vt-wrt-P}, for the given $\epsilon$ there exists a constant $\delta_0>0$ such that 
\bgeq
\qquad \dd_{\inmat{Prok}}\left(
\mathbbm{1}_{\inf_{f\in \F} R_P(f)},\mathbbm{1}_{\inf_{f\in \F} R_Q(f)}\right) \leq \frac{\epsilon}{3}, \forall 
Q \in {\cal M} \subset  {\cal M}_{Z,\kappa}^{\phi^p} \; \inmat{with}  \; \dd_{\phi^p}(P,Q) \leq \delta_0.
\label{eq:stat-robust-vt_n-SP-2}
\edeq
So we are left to show that
\bgeqn
\dd_{\inmat{Prok}}\left(P^{\otimes N}\circ \hat{\vt}_N(\cdot,\lambda_N)^{-1},
\mathbbm{1}_{\inf_{f\in \F} R_P(f)}\right)
\leq \frac{\epsilon}{3}
\label{eq:stat-robust-vt_n-SP-3}
\edeqn
and
\bgeqn
\dd_{\inmat{Prok}}\left(
\mathbbm{1}_{\inf_{f\in \F} R_Q(f)},
Q^{\otimes N}\circ \hat{\vt}_N(\cdot,\lambda_N)^{-1}
\right)\leq \frac{\epsilon}{3}
\label{eq:stat-robust-vt_n-SP-4}
\edeqn
for $N$ sufficiently large. By Strassen's theorem \cite{Hub81},
(\ref{eq:stat-robust-vt_n-SP-3}) and (\ref{eq:stat-robust-vt_n-SP-4})
are implied respectively by
\bgeqn
P^{\otimes N}
\left[
\vec{z}^N: \left|
\hat{\vt}_N(\vec{z}^N,\lambda_N) -\inf_{f\in \F} R_P(f)
\right|
\geq \frac{\epsilon}{3}
\right]
\leq \frac{\epsilon}{3}
\label{eq:stat-robust-vt_n-SP-5}
\edeqn
and
\bgeqn
Q^{\otimes N}
\left[
\tilde{\vec{z}}^N: \left|
\hat{\vt}_N(\tilde{\vec{z}}^N,\lambda_N) -\inf_{f\in \F} R_Q(f)
\right|
\geq \frac{\epsilon}{3}
\right]
\leq \frac{\epsilon}{3}.
\label{eq:stat-robust-vt_n-SP-6}
\edeqn
Using the definition of the optimal values,
 (\ref{eq:stat-robust-vt_n-SP-5}) and (\ref{eq:stat-robust-vt_n-SP-6}) can be 
rewritten respectively  as 
\bgeqn
P^{\otimes N}
\left[
\vec{z}^N: 
\left|
\inf_{f\in\F}  \bbe_{P_N}\{\left[ c(z,f(x))\right]
+ \lambda_N \left\|f\right\|_k^2\}
-\inf_{f\in \F} R_P(f)
\right|
\geq \frac{\epsilon}{3}
\right]
\leq \frac{\epsilon}{3}
\label{eq:stat-robust-vt_n-SP-7}
\edeqn
and
\bgeqn
Q^{\otimes N}
\left[
\tilde{\vec{z}}^N: \left|
\inf_{f\in\F} \bbe_{Q_N}\{\left[c(z,f(x))\right]
+ \lambda_N \left\|f\right\|_k^2\}
-\inf_{f\in \F} R_Q(f)
\right|
\geq \frac{\epsilon}{3}
\right]
\leq \frac{\epsilon}{3}.
\label{eq:stat-robust-vt_n-SP-8}
\edeqn
Note that we may set 
$N_0\in\mathbb{N}$  sufficiently large such that
$
\lambda_N\leq \frac{\epsilon}{6\beta^2}
$
for all $N\geq N_0$.
Consequently the two inequalities above are implied by
\bgeqn
P^{\otimes N}
\left[
\vec{z}^N: 
\left|
\inf_{f\in\F} R_{P_N}(f)
-\inf_{f\in \F} R_P(f)
\right|
\geq \frac{\epsilon}{6}
\right]
\leq \frac{\epsilon}{3}
\label{eq:stat-robust-vt_n-SP-9}
\edeqn
and
\bgeqn
Q^{\otimes N}
\left[
\tilde{\vec{z}}^N: 
\left|
\inf_{f\in\F} R_{Q_N}(f)
-\inf_{f\in \F} R_Q(f)
\right|
\geq \frac{\epsilon}{6}
\right]
\leq \frac{\epsilon}{3},
\label{eq:stat-robust-vt_n-SP-10}
\edeqn
or equivalently
\bgeqn
P^{\otimes N}
\left[
\vec{z}^N: 
\left|
\hat{\vt}_N(\vec{z}^N)-\vt(P)
\right|
\geq \frac{\epsilon}{6}
\right]
\leq \frac{\epsilon}{3}
\label{eq:stat-robust-vt_n-SP-9-vt}
\edeqn
and
\bgeqn
Q^{\otimes N}
\left[
\tilde{\vec{z}}^N: 
\left|
\hat{\vt}_N(\tilde{\vec{z}}^N)-\vt(Q)
\right|
\geq \frac{\epsilon}{6}
\right]
\leq \frac{\epsilon}{3}.
\label{eq:stat-robust-vt_n-SP-10-vt}
\edeqn
By Theorem \ref{t-cont-vt-wrt-P}, there exists a constant $\delta>0$ such that when
$\dd_{\phi^p}(P',P)< 2\delta$, 
$\left|
\vt(P')-\vt(P)
\right|< \frac{\epsilon}{12}$.
On the other hand, 
it follows by \cite[Corollary 3.5]{KSZ12} that
$(\M^{\phi^p}_{Z,\kappa},\dd_{\phi})$ has the UGC property which implies that
\bgeqn
Q^{\otimes N}\left[ \dd_{\phi^p}(Q_N,Q)\geq \delta\right] \leq \frac{\epsilon}{3}
\label{UGC-prok-P-N-x=-to-Px}
\edeqn
for all $Q\in \M^{\phi^p}_{Z,\kappa}$ including $Q=P$. This shows (\ref{eq:stat-robust-vt_n-SP-9-vt})
when $N_0$ is chosen sufficiently large.
To show (\ref{eq:stat-robust-vt_n-SP-10-vt}), let
$\dd_{\phi^p}(Q,P)\leq \delta$. Then
\bgeqn
\frac{\epsilon}{3} 
&\geq& Q^{\otimes N}
\left[\tilde{\vec{z}}^N:  \dd_{\phi^p}(Q_N,Q)
\geq \delta \right]\\
&\geq& 
Q^{\otimes N}
\left[
\tilde{\vec{z}}^N: \dd_{\phi^p}(Q_N,P) \geq \delta+\dd_{\phi^p}(Q,P)\right]\nonumber\\
&\geq& 
Q^{\otimes N}
\left[\tilde{\vec{z}}^N: 
\dd_{\phi^p}(Q_N,P) 
\geq 2\delta \right]\nonumber\\
&\geq& 
Q^{\otimes N}
\left[\tilde{\vec{z}}^N
: 
\left|
\vt(Q_N)-\vt(P)
\right|
\geq \frac{\epsilon}{12}
\right]\nonumber\\
&\geq& 
Q^{\otimes N}
\left[
\tilde{\vec{z}}^N: 
\left|
\vt(Q_N)-\vt(Q)
\right|
\geq|\vt(P)-\vt(Q)|
+\frac{\epsilon}{12}
\right]\nonumber\\
&\geq& Q^{\otimes N}
\left[
\tilde{\vec{z}}^N: 
\left|
\vt(Q_N)-\vt(Q)
\right|
\geq \frac{\epsilon}{6}
\right], \forall Q\in \M^{\phi^p}_{Z,\kappa}.\nonumber
\label{eq:Uniform-consist}
\edeqn
The conclusion follows.
\hfill $\Box$

We make a few comments about the conditions and results of this theorem.

First, the set
 $ \M_{Z,\kappa}^{\phi^p}$
differs from
 ${\cal M}_{Z}^{\phi^p}$
in that the former imposes a bound for the moment value uniformly for all $P\in \M_{Z,\kappa}^{\phi^p}$ whereas
the latter does not have such
uniformity. 
This is because we need
the UGC property of $(\M_{Z,\kappa}^{\phi^p},\dd_{\phi})$ in order for us to apply  
\cite[Corollary 3.5]{KSZ12}.
For example, 
in the least squares regression model with polynomial kernel, we have 
$$
\M_{Z,\kappa}^{\phi^p} = \left\{P\in \mathscr{P}(Z): \int_Z \left[\|y\|^2+\beta^2(\gamma \|x\|^2+1)^{d}\right]^pP(dz)<\kappa \right\}.
$$
In the case of Gaussian kernel or Laplacian kernel,
$$
\M_{Z,\kappa}^{\phi^p} = \left\{P\in \mathscr{P}(Z): \int_Z \|y\|^{2p} P(dz)<\kappa \right\}.
$$

Second, by  (\ref{eq:stat-robust-vt_n-SP-10-vt}),
we can obtain for any $\epsilon>0$, there exist constants $\delta>0$ and $N_0 \in \mathbb{N}$ such that
\bgeq
Q \in \mathcal{M}, \dd_{\phi^p}(P,Q) \leq \delta
\Longrightarrow
Q^{\otimes N}\left[
\tilde{\vec{z}}^N: \left|
\vt(Q)-\vt(Q_N)
\right|\geq \frac{\epsilon}{6}
\right] \leq \frac{\epsilon}{3}
\edeq
for $N \geq N_0$.
This implies  uniform convergence of $\vt(Q_N)$ to $\vt(Q)$ for
all $Q$ near $P$ as opposed to pointwise convergence (for each fixed $Q$) in stochastic programming.
The uniformity
does not come out 
for free: it restricts both $P$ and $Q$
to the $\phi$-weak topological space of probability measures. 


Third, in practice, since $P$ is unknown, it is difficult to identify $\delta$ for a specified $\epsilon$.
The usefulness of  (\ref{eq:stat-robust-vt_n-SP}) should be understood as that it provides a theoretical guarantee:
if  the training data are generated by some probability distribution $Q$ which is close to
the true distribution $P$, and $Q$ satisfies moment condition (\ref{eq:M-k-kappa}) (which may be examined
through empirical data), then the optimal value obtained with the perceived data is close to
the one with real data.
There are potentially two ways to move forward the research.
One is to derive quantitative statistical robustness
under some additional conditions in which case the relationship between $\epsilon$ and $\delta$ may be explicitly established, we will come back to this in 
the next subsection.
The other is to use the training data to construct an ambiguity set of
probability distributions and use the latter to develop a model which is robust both
in preference and in brief. This will effectively create
a robust mechanism to mitigate the risk arising from noise in
perceived data. We leave this for future research.

\subsection{Quantitative robustness}

In the previous section, there is no explicit relationship between $\epsilon$ and $\delta$ in the 
qualitative robustness result. In this section, we address the issue
under the following additional conditions.

\begin{assu}
\label{A:cost-Lip}
The cost function $c(z,f(x))$ satisfies the following property:
\bgeqn
\label{eq:local-Lip-assu-0} 
|c(z,f(x)) - c(z',f(x'))| 
\leq c_p(z,z')\|z-z'\|, \forall z,z'\in Z, 
f \in \F,
\edeqn
where $c_p(z,z') :=\max\{1,\|z\|,\|z'\|\}^{p-1}$ and $p\geq 1$ is a fixed positive number.
\end{assu}

To see how the assumption may be satisfied,
we consider the case that $c(z,f(x))$ 
is locally Lipschitz continuous with modulus being bounded by $L(z)$, then
\bgeq
|c(z,f(x)) - c(z',f(x'))| \leq \max\{L(z),L(z')\}(\|z-z'\|+|f(x)-f(x')|), \forall z,z'\in Z.
\edeq
Under Assumption \ref{A:kernel} (b) and the calmness condition
in Remark \ref{r:calm},
\bgeq
|f(x)-f(x')| = |\langle f, k(\cdot,x)\rangle - \langle f, k(\cdot,x')\rangle|
\leq  \beta \|k(\cdot,x)-k(\cdot,x')\|_k \leq \beta g(\|x-x'\|).
\edeq
Consequently we have
\bgeqn
|c(z,f(x)) - c(z',f(x'))| \leq \max\{L(z),L(z')\}(\|z-z'\|+\beta g(\|x-x'\|)), 
\forall z,z'\in Z.\nonumber\\
\label{eq:local-Lip-assu-1} 
\edeqn
In Example \ref{ex:cost-quant-local-Lip}, we will 
explain in detail how $L(\cdot)$ may be figured out and in a combination with specific form of 
function $g(\cdot)$, inequality (\ref{eq:local-Lip-assu-1}) leads to inequality (\ref{eq:local-Lip-assu-0}) 
for some specific cost functions and kernel functions in regression models.

We now return to our discussion on the quantitative description 
of the discrepancy between $P^{\otimes N}\circ \hat{\vt}_N(\cdot,\lambda_N)^{-1}$
and $Q^{\otimes N}\circ \hat{\vt}_N(\cdot,\lambda_N)^{-1}$. Our idea 
is to use Kantorovich metric to measure the difference, i.e.,
$
\dd_{K,1} \left(P^{\otimes N}\circ \hat{\vt}_N(\cdot,\lambda_N)^{-1},
Q^{\otimes N}\circ \hat{\vt}_N(\cdot,\lambda_N)^{-1}\right)
$, which can be converted to 
the estimate of the difference  
between $P^{\otimes N}$ and $Q^{\otimes N}$ under some metric 
by the $\zeta$-metric of $P$ and $Q$. The next technical 
result prepares for such a conversion. 

\begin{lema}
\label{lem:FM-metric-estimate}
Let $\vec{z}:=(z^1,\cdots,z^N) \in Z^{\otimes N}$ and
$$
\Psi:=\left\{\psi: Z^{\otimes N} \to \R:
|\psi(\vec{\tilde{z}})-\psi(\vec{\hat{z}})|
\leq \frac{1}{N}\sum_{j=1}^N c_{p}(\tilde{z}^j,\hat{z}^j)\|\tilde{z}^j-\hat{z}^j\|
\right\}.
$$
Let 
$
\dd_{\Psi} (P^{\otimes N},Q^{\otimes N})= \sup_{\psi\in \Psi} \left|\int_{Z}\psi(z)P^{\otimes N}(dz)-\int_{Z}\psi(z)Q^{\otimes N}(dz)\right|.
$
Then
$$
\dd_{\Psi} (P^{\otimes N},Q^{\otimes N}) \leq \zeta_{p} (P,Q).
$$
\end{lema}


\noindent\textbf{Proof.}
The result is established in \cite[Lemma 4.1]{WXM2020} which is an extension of \cite[Lemma 1]{GuX18}
(which is presented when $p=1$).
Here we include a proof for self-containedness.
Let 
$\vec{z}^j:=\{z^1,\cdots,z^j\}$ and
$\vec{z}^{-j}:=\{z^1,\cdots,z^{j-1},z^{j+1},\cdots,z^N\}$.
For any $P_1,\cdots,P_N \in \mathscr{P}(Z)$ and any $j\in \{1,\cdots,N\}$,
denote
$$
P_{-j}(d \vec{z}^{-j}):=P_1(d z^1)\cdots P_{j-1}(d z^{j-1})P_{j+1}(d z^{j+1})\cdots P_N(d z^N)
$$
and
$h_{\vec{z}^{-j}}({z}^j):=\int_{Z^{\otimes (N-1)}}
\psi(\vec{z}^{-j},z^j) P_{-j}(d \vec{z}^{-j})$.
Then
\bgeq
|h_{\vec{z}^{-j}}(\tilde{z}^j)-h_{\vec{z}^{-j}}(\hat{z}^j)|
&\leq& \int_{Z^{\otimes (N-1)}}
\left|\psi(\vec{z}^{-j},\tilde{z}^j)-\psi(\vec{z}^{-j},\hat{z}^j)\right|
P_{-j}(d \vec{z}^{-j})\\
&\leq& \int_{Z^{\otimes (N-1)}}
\frac{1}{N}c_{p}(\tilde{z}^j,\hat{z}^j)\|\tilde{z}^j-\hat{z}^j\|P_{-j}(d \vec{z}^{-j})\\
&\leq& \frac{1}{N}c_{p}(\tilde{z}^j,\hat{z}^j)\|\tilde{z}^j-\hat{z}^j\|.
\edeq
Let $\mathcal{W}$ denote the set of functions $h_{\vec{z}^{-j}}({z}^j)$ generated by $\psi\in \Psi$.
By the definition of $\dd_{\Psi}$ and the $p$-th order Fortet-Mourier metric,
\bgeqn
\dd_{\Psi} (P_{-j} \times \tilde{P}_j,P_{-j} \times \hat{P}_j)
&=&\sup_{\psi \in \Psi}
\left|
\int_{Z}\int_{Z^{\otimes (N-1)}}
\psi(\vec{z}^{-j},z^j) P_{-j}(d \vec{z}^{-j}) \tilde{P}_j(d z^j)\right.\nonumber\\
&&-\left.
\int_{Z}\int_{Z^{\otimes (N-1)}}
\psi(\vec{z}^{-j},z^j) P_{-j}(d \vec{z}^{-j}) \hat{P}_j(d z^j)\right|\nonumber\\
&=& \sup_{h_{\vec{z}^{-j}} \in {\cal W}} \left|
\int_{Z}h_{\vec{z}^{-j}}({z}^j)\tilde{P}_j(d z^j)
-\int_{Z}h_{\vec{z}^{-j}}({z}^j)\hat{P}_j(d z^j) \right|\nonumber\\
&\leq& \frac{1}{N}\zeta_{p}(\tilde{P}_j,\hat{P}_j),
\label{eq:productP_j-tildeP_j}
\edeqn
where the inequality is due to $N h_{\vec{z}^{-j}}(z^j)\in \mathcal{F}_{p}(Z)$ and the definition of $\zeta_{p}(P,Q)$. 
Finally, by the triangle inequality of the pseudo-metric, we have
\bgeq
\dd_{\Psi} \left(P^{\otimes N},Q^{\otimes N}\right)
&\leq&\dd_{\Psi} \left(P^{\otimes N}, P^{\otimes (N-1)}\times Q\right)
+\dd_{\Psi} \left(P^{\otimes (N-1)}\times Q, P^{\otimes (N-2)}\times Q^{\otimes 2}\right)\\
&&+ \cdots+ \dd_{\Psi} \left(P\times Q^{\otimes (N-1)}, Q^{\otimes N} \right)\\
&\leq& \frac{1}{N} \zeta_{p}(P,Q) \times N= \zeta_{p}(P,Q).
\edeq
The proof is complete.
\hfill $\Box$

With Lemma \ref{lem:FM-metric-estimate}, we are ready to state our main result.

\begin{thm}[Quantitative statistical robustness]
\label{T-Quant-SR}
Let $\phi(z)$ be defined as in Assumption~\ref{A:cost-1}
and
$
{\cal M}_Z^\phi = \left\{P'\in \mathscr{P}(Z): \int_{Z} \phi(z)P'(dz)<\infty\right\}.
$
Under Assumptions~\ref{A:kernel} (b), \ref{A:cost-1} (a) and \ref{A:cost-Lip},
\bgeqn
 \dd_{K,1} \left(P^{\otimes N}\circ \hat{\vt}_N(\cdot,\lambda_N)^{-1},
Q^{\otimes N}\circ \hat{\vt}_N(\cdot,\lambda_N)^{-1}\right) \leq  \zeta_p (P,Q)
\label{eq:quant-SR}
\edeqn
 for any $N \in \mathbb{N}$ and any $P, Q\in {\cal M}_Z^{\phi}$,
 where $p$ is defined as in Assumption \ref{A:cost-Lip}.
 In the case when $p=1$,
\bgeqn
 \dd_{K,1} \left(P^{\otimes N}\circ 
\hat{\vt}_N(\cdot,\lambda_N)^{-1},
Q^{\otimes N}\circ \hat{\vt}_N(\cdot,\lambda_N)^{-1}\right) \leq  \dd_{K,Z} (P,Q).
\label{eq:quant-SR-Kant}
\edeqn
\end{thm}

\noindent\textbf{Proof.}
By definition
\bgeqn
\label{eq:stat-robust-vt_n-thm-pseudo}
&& \dd_{K,1} \left(P^{\otimes N}\circ \hat{\vt}_N(\cdot,\lambda_N)^{-1},
Q^{\otimes N}\circ \hat{\vt}_N(\cdot,\lambda_N)^{-1}\right) \\
&=&\sup_{g \in \mathscr{G}}\left|
\int_\R g(t) P^{\otimes N}\circ \hat{\vt}_N(\cdot,\lambda_N)^{-1}(dt) -
\int_\R g(t) Q^{\otimes N}\circ \hat{\vt}_N(\cdot,\lambda_N)^{-1}(dt)
\right|\nonumber\\
&=&
\sup_{g \in \mathscr{G}}\left|
\int_{Z^{\otimes N}}g(\hat{\vt}_N(
\vec{z}^N,\lambda_N)) P^{\otimes N}
(d\vec{z}^N) -
\int_{Z^{\otimes N}} g(\hat{\vt}_N(
\vec{z}^N,\lambda_N)) Q^{\otimes N}
(d\vec{z}^N)
\right|,\nonumber
\edeqn
where we write $\vec{z}^N$ for $(z^1,\cdots,z^N)$
and
$\hat{\vt}(\vec{z}^N,\lambda_N)$ for
$\hat{\vt}_N$ to indicate its dependence on $z^1,\cdots,z^N$.
To see the well-definiteness of the pseudo-metric, we note that
for each $g\in \mathscr{G}$,
\bgeqn\label{eq:quanti-well-defi1}
|g(\hat{\vt}_N(\vec{z}^N,\lambda_N))|
\leq |g(\hat{\vt}_N(\vec{z}_0^N,\lambda_N))|
+ |\hat{\vt}_N(\vec{z}^N,\lambda_N)-\hat{\vt}_N(\vec{z}_0^N,\lambda_N)|,
\edeqn
where $\vec{z}^N_0 \in Z^{\otimes N}$ is fixed.
By the definition of $\hat{\vt}(\vec{z}^N,\lambda_N)$, we have
\bgeq
|\hat{\vt}_N(\vec{z}^N,\lambda_N)|
&=&\left|\displaystyle{
\min_{f\in \F}   }  \;\;
\frac{1}{N}\sum_{j=1}^N \left(c(z^j,f(x^j))+\lambda_N \|f\|_k^2\right)\right|
\leq
\frac{1}{N}\sum_{j=1}^N \phi(z^j)+\lambda_N\beta^2.
\edeq
Thus
\bgeqn\label{eq:quanti-well-defi2}
\int_{Z^{\otimes N}} |\hat{\vt}_N(\vec{z}^N,\lambda_N)|
P^{\otimes N}
(d\vec{z}^N)
&\leq&
\int_{Z^{\otimes N}}
\frac{1}{N}\sum_{j=1}^N \phi(z^j)P^{\otimes N}
(d\vec{z}^N)+\lambda_N\beta^2 \nonumber\\
&=&\int_{Z}\phi(z)P(dz)+\lambda_N\beta^2<\infty, \forall P\in {\cal M}_{Z}^{\phi},
\edeqn
where the
equality holds due to
the fact that $z^1,\cdots,z^N$ are i.i.d..
The same inequality can be established for $\int_{Z^{\otimes N}} |\hat{\vt}_N(\vec{z}_0^N,\lambda_N)|
P^{\otimes N}
(d\vec{z}^N)$.
Combining (\ref{eq:quanti-well-defi1}) and (\ref{eq:quanti-well-defi2}),
we deduce that
$$
\int_{Z^{\otimes N}}g(\hat{\vt}_N(\vec{z}^N,\lambda_N)) P^{\otimes N}
(d\vec{z}^N)<\infty, \forall P \in {\cal M}_{Z}^{\phi}.
$$
The same argument can be made
on $\int_{Z^{\otimes N}}
g(\hat{\vt}_N(\vec{z}^N,\lambda_N)) Q^{\otimes N}
(d\vec{z}^N)$ for $Q \in {\cal M}_{Z}^{\phi}$.

Next, we show (\ref{eq:quant-SR}).
We do so by applying Lemma \ref{lem:FM-metric-estimate} to the right hand side of
(\ref{eq:stat-robust-vt_n-thm-pseudo}).
To this end, we need to verify the condition of the lemma.
Define  $\psi: Z^{\otimes N} \to \R$ by
$
\psi(\vec{z}^N):=g(\hat{v}(\vec{z}^N,\lambda_N)). 
$
Since 
$g$ is Lipschitz continuous with modulus bounded by $1$,
we have
\bgeq
&&\left|\psi(\tilde{\vec{z}}^N)-\psi(\hat{\vec{z}}^N)
\right|\\
&=&
|g(\hat{\vt}_N(\vec{\tilde{z}}^N,\lambda_N))
-g(\hat{\vt}_N(\vec{\hat{z}}^N,\lambda_N))|\\
&\leq& 
|\hat{\vt}_N(\tilde{\vec{z}}^N,\lambda_N)
-\hat{\vt}_N(\hat{\vec{z}}^N,\lambda_N)| \\
&=& 
\left|
\min_{f\in \F}     \;\;
\frac{1}{N}\sum_{j=1}^N \left(c(\tilde{z}^j,f(x^j))+\lambda_N \|f\|_k^2\right)
-
\min_{f\in \F}   \;\;
\frac{1}{N}\sum_{j=1}^N \left(c(\hat{z}^j,f(x^j))+\lambda_N \|f\|_k^2\right)
\right|\\
&\leq & \frac{1}{N} \sum_{j=1}^N\sup_{f\in \F}  |c(\tilde{z}^j,f(\tilde{x}^j))-c(\hat{z}^j,f(\hat{x}^j))|\\
&\leq&\frac{1}{N} \sum_{j=1}^N c_p(\tilde{z}^j,\hat{z}^j)\|\tilde{z}^j-\hat{z}^j\|,
\edeq
which means that $\psi$ is in the set of functions $\Psi$ in  Lemma \ref{lem:FM-metric-estimate}.
The rest follows from application of the lemma to (\ref{eq:stat-robust-vt_n-thm-pseudo}).
\hfill $\Box$

The strength  of 
 Theorem \ref{T-Quant-SR}
lies in the fact that it gives rise to 
an explicit quantitative relationship between $ \dd_{K,1} \left(P^{\otimes N}\circ \hat{\vt}_N(\cdot,\lambda_N)^{-1},
Q^{\otimes N}\circ \hat{\vt}_N(\cdot,\lambda_N)^{-1}\right)$ and
$\zeta_p (P,Q)$. This is benefited partially from use of the dual representation of the Kantorovich metric
in the quantification of the former and partially from use of Fortet-Mourier metric for quantification of the latter.
As noted immediately after Definition \ref{D-Fort-Mou-metric}, $\zeta_p(P,Q)$ may be estimated 
via sample data, which means the 
error bound established in (\ref{eq:quant-SR}) is practically obtainable and this is a significant
step forward from the qualitative robustness result. 
Note also that  both $\dd_\phi$ and
$\zeta_p$ capture (restrict) the tail behaviour of $P$ but there is no explicit relationship between the two metrics as far as we are concerned: the former provides weak convergence of {\em each}  measurable function dominated by $\phi$ whereas the latter
requires {\em uniform} convergence of a class of locally Lipschitz continuous functions with specified 
rate of growth.
Finally, we note that the error bound does not depend on the regularization parameters because from the proof we can see that
the regularization terms are cancelled. It does not mean that the parameter has no effect on the statistical performance of
the empirical risk, rather it means the error bound does not capture such effect.

The next example illustrates how the theorem works in some concrete regression models. 

\begin{ex}
\label{ex:cost-quant-local-Lip}
Consider 
the least squares regression model, where
$c(z,f(x))=\frac{1}{2}|y-f(x)|^2$. We have 
\bgeq
|c(z,f(x))-c(z',f(x'))|
&=&\frac{1}{2}\left| |y-f(x)|^2-|y'-f(x')|^2 \right|\\
&\leq&\frac{1}{2}\left(|y|+|f(x)|+|y'|+|f(x')|)
(|y-y'|+|f(x)-f(x')|\right).
\edeq
Under Assumption \ref{A:kernel} (b) and the calmness condition in Remark \ref{r:calm},
\bgeq
|f(x)| \leq \|f\|_k \|k(x,\cdot)\|_k \leq \beta \|k(x,\cdot)\|_k
=\beta \sqrt{k(x,x)}, \forall f\in \F
\edeq
and
\bgeq
|f(x)-f(x')| & =& |\langle f, k(\cdot,x)\rangle - \langle f, k(\cdot,x')\rangle|
\leq  \beta \|k(\cdot,x)-k(\cdot,x')\|_k
\leq \beta g(\|x-x'\|).
\edeq
Let $\eta(z):=|y|+\beta \sqrt{k(x,x)}$.
Then,
\bgeq
\label{eq:C-least-sqr}
|c(z,f(x))-c(z',f(x'))|
\leq 
\max\left\{\eta(z),\eta(z')\right\}
(|y-y'|+\beta g(\|x-x'\|)).
\edeq
\begin{itemize}

\item In the case of linear kernel,
$\eta(z)=|y|+\beta \|x\| 
\leq \max\{1,\beta\}\|z\|$, $g(t)=t$, and
\bgeq
|c(z,f(x))-c(z',f(x'))|
\leq 
\max\{1,\beta\}^2 \max\left\{1,\|z\|,\|z'\|\right\}
\|z-z'\|.
\edeq
By Theorem \ref{T-Quant-SR}, 
$
 \dd_{K,1} \left(P^{\otimes N}\circ \hat{\vt}_N^{-1},
Q^{\otimes N}\circ \hat{\vt}_N^{-1}\right) \leq 
\max\{1,\beta\}^2 \zeta_2 (P,Q)$
 for all $N \in \mathbb{N}$ and any $P, Q\in {\cal M}_Z^{\phi}$, where
  $\phi(z)=\|y\|^2+\beta^2\|x\|^2$.

\item 
In the case of Gaussian kernel, $\eta(z)=|y| 
\leq \|z\|$, $g(t)=\max\{\sqrt{2\gamma},1\}t$, and
\bgeq
|c(z,f(x))-c(z',f(x'))
\leq 
\max\{\sqrt{2\gamma},1\} \max\left\{1,\|z\|,\|z'\|\right\}
\|z-z'\|.
\edeq
By Theorem \ref{T-Quant-SR},
$ 
 \dd_{K,1} \left(P^{\otimes N}\circ \hat{\vt}_N^{-1},
Q^{\otimes N}\circ \hat{\vt}_N^{-1}\right) \leq 
\max\{\sqrt{2\gamma},1\} \zeta_2 (P,Q)
$ 
 for all $N \in \mathbb{N}$ and any $P, Q\in {\cal M}_Z^{\phi}$, where  $\phi(z)=\|y\|^2$.

\item In the case of polynomial kernel,
$\eta(z)=|y|+\beta \sqrt{(\gamma \|x\|^2+1)^{d}}$.
For fixed $z$ and $z'$, let $R:=\max\{1,\|z\|,\|z'\|\}$.
Then
\bgeq
&&\|k(\cdot,x)-k(\cdot,x')|_k\\
&\leq&
\left\{
\begin{array}{ll}
\max\{\frac{1}{2R}\sqrt{2(\gamma R^2 +1)^d},1\}\|x-x'\|, & \inmat{if} \; d\, \inmat{is even},\\
\max\{\frac{1}{2R}\sqrt{2(\gamma R^2 +1)^d}-
2(1-\gamma R^2)^d,1\}\|x-x'\|, & \inmat{if} \; d\, \inmat{is odd},
\end{array}
\right.\\
&\leq& 
\left\{
\begin{array}{ll}
\max\{(\frac{1+\gamma}{2})^{d/2},1\} \max\{1,\|z\|,\|z'\|\}^{d-1}\|x-x'\|, & \inmat{if} \; d\, \inmat{is even},\\
\max\left\{
2(\frac{1+\gamma}{2})^{d/2},4,
4\gamma^d\right\}
\max\{1,\|z\|,\|z'\|\}^{2d}
\|x-x'\|, & \inmat{if} \; d\, \inmat{is odd}.
\end{array}
\right.
\edeq
The last inequality is due to the fact that 
$a-b \leq \max\{2a,-2b\}$ for any two numbers 
$a,b$ where  $a>0$ and $b$
could be either negative or positive. 
Let
\bgeq
A_1&:=&(1+\beta(\gamma+1)^{d/2})
\max\left\{\beta(\frac{1+\gamma}{2})^{d/2},\beta,1\right\},\\
A_2&:=&(1+\beta(\gamma+1)^{d/2})
\max\left\{2\beta(\frac{1+\gamma}{2})^{d/2},4\beta,
4\beta\gamma^d,1\right\}.
\edeq
Then
\bgeq
|c(z,f(x))-c(z',f(x'))|
\leq 
\left\{
\begin{array}{ll}
A_1\max\left\{1,\|z\|,\|z'\|\right\}^{2d-1}
\|z-z'\|,  & \inmat{if} \; d\, \inmat{is even},\\
A_2\max\left\{1,\|z\|,\|z'\|\right\}^{3d}
\|z-z'\|, &\inmat{if} \; d\, \inmat{is odd}.
\end{array}
\right.
\edeq
By Theorem \ref{T-Quant-SR}
\bgeq
 \dd_{K,1} \left(P^{\otimes N}\circ \hat{\vt}_N^{-1},
Q^{\otimes N}\circ \hat{\vt}_N^{-1}\right) \leq 
\left\{
\begin{array}{ll}
A_1\zeta_{2d} (P,Q), & \inmat{if} \; d\, \inmat{is even},\\
A_2 \zeta_{3d+1} (P,Q), &\inmat{if} \; d\, \inmat{is odd},
\end{array}
\right.
\label{eq:quant-SR-3}
\edeq
 for all $N \in \mathbb{N}$ and any $P, Q\in {\cal M}_Z^{\phi}$, where $\phi(z)=\|y\|^2+\beta^2(\gamma \|x\|^2+1)^{d}$.

\end{itemize}

We can derive similar results for 
the regression models with  $\epsilon$-insensitive loss function
$c(z,f(x))=\max\{0,|y-f(x)|-\epsilon\}$,
hinge loss $c(z,f(x))=\max\{0,1-(y-f(x))\}$,
and log-loss function
$c(z,f(x))=\log(1+e^{-(y-f(x))})$ respectively, we omit the details.
\end{ex}

\begin{rem}
It might be interesting to 
study the discrepancy between 
$f_N^{\lambda_N}(P_N)$ and $f_N^{\lambda_N}(Q_N)$.
To this end, 
we assume that $c(z,f(x))$ is strong convex in $f$ 
for almost all $z$. In such a case, 
$R(f)=\bbe_P[c(z,f(x))]$ is also strongly convex and so is 
$R(f)+\lambda\|f\|_k$,
which
implies that  problem (\ref{eq:ML-Rf-min})
and the regularized problem (\ref{eq:ML-saa-r}) have a unique solution.
Moreover, the strong convexity 
implies that  problem (\ref{eq:ML-saa-r})  satisfies second order growth condition at $f_N^{\lambda_N}(P_N)$, 
that is,
there exists a positive constant $\alpha$
such that
\bgeq
R_{P_N}^{\lambda_N}(f)- \vt(P_N,\lambda_N)
\geq \alpha \|f-f_N^{\lambda_N}(P_N)\|_k^2, \forall f\in \F.
\edeq
By virtue of \cite[Lemma 3.8]{LX13}, we can use the inequality to obtain 
\bgeq
\|f_N^{\lambda_N}(P_N)-f_N^{\lambda_N}(Q_N)\|_k 
&\leq& \sqrt{\frac{3}{\alpha}\sup_{f\in \F}|\bbe_{P_N}[c(z,f(x))]
-\bbe_{Q_N}[c(z,f(x))]|}\\
&\leq&\sqrt{\frac{3}{\alpha}\bbe_{P_N\times Q_N}[c_p(\hat{z},\tilde{z})\|\hat{z}-\tilde{z}\|]}.
\edeq
Since $\bbe_{P_N\times Q_N}[c_p(\hat{z},\tilde{z})\|\hat{z}-\tilde{z}\|]
-
\bbe_{P\times Q}[c_p(\hat{z},\tilde{z})\|\hat{z}-\tilde{z}\|]\to 0$ as $Q\to P$ and $N$ goes to infinity, 
then 
$$
\|f_N^{\lambda_N}(P_N)-f_N^{\lambda_N}(Q_N)\| \to 0.
$$
However, we are unable to establish the kind of estimation in (\ref{eq:quant-SR}) for the optimal solutions 
because of the non-linearity of the bound $$\sqrt{\frac{3}{\alpha}\sup_{f\in \F}|\bbe_{P_N}[c(z,f(x))]
-\bbe_{Q_N}[c(z,f(x))]|}$$ for $\|f_N^{\lambda_N}(P_N)-f_N^{\lambda_N}(Q_N)\|_k $
in terms of the difference of the function values.


\end{rem}

\section{Uniform consistency}

In this section, we move on to investigate convergence of 
$\vt(P_N,\lambda_N)$ to $\vt(P)$ as $N\to \infty$ and $\lambda_N\to 0$.
We proceed the investigation in two steps: first pointwise convergence, i.e., for each fixed $P\in \mathscr{P}(Z)$ and 
then uniform convergence for all $P$ over a subset ${\cal M}$ of $\mathscr{P}(Z)$. To this end, we introduce the following
assumption on the cost function.


\begin{assu} 
\label{A:cost-Hold}
There exist a measurable function $r(\cdot): Z \rightarrow \R_+$ and a constant $\nu\in (0,1]$ such that
\bgeqn
|c(z,f(x))-c(z,g(x))| \leq r(z) \|f-g\|_\infty^\nu,  \forall f,g\in \F, z \in Z.
\label{eq:Holder-condition-c}
\edeqn
\end{assu}

The assumption requires $c(z,\cdot)$ to be H\"older continuous over ${\cal F}$ uniformly
for $z\in Z$. It should be distinguished from Assumption \ref{A:cost-Lip} which requires
$c(z,f(x))$ to be locally Lipschitz continuous in $z$ for all $f\in \F$.
The assumption is satisfied by all of the loss functions in regression models that we listed at the beginning of
Section 2.

\begin{thm}[Consistency of $\vt(P_N,\lambda_N)$]
\label{thm-vt-consist-ML}
Let Assumptions \ref{A:kernel}, \ref{A:cost-1}  and \ref{A:cost-Hold} hold.
 Then for any  $\delta>0$, there
exist  positive constants $\e<\delta/6$,
$\alpha(\e,\delta)$ and  $\gamma(\e,\delta)$, independent
of $N$ and a positive number $N_0$ such that
 \bgeqn
\label{eq:lemma-saa-approxi-1}
\quad P^{\otimes N} \left(
\sup_{f\in {\cal F}}
|\bbe_{P_N}[c(z,f(x))]+\lambda_N \|f\|_k^2 - \bbe_{P}[c(z,f(x))|\geq \delta
\right ) \leq
\alpha(\e,\delta) e^{-N\gamma (\e,\delta)}
\edeqn
when $N\geq N_0$ and $\lambda_N \leq \epsilon/\beta^2$
and hence
 \bgeqn
\label{eq:lemma-saa-approxi-2}
P^{\otimes N} \left(
|\vt(P_N,\lambda_N) - \vt(P) |\geq \delta
\right) \leq
\alpha(\e,\delta) e^{-N\gamma(\e,\delta)}
\edeqn
and
 \bgeqn
\label{eq:lemma-saa-approxi-2-a}
P^{\otimes N} \left(
| \bbe_P[c(z,f_N^{\lambda_N}(x))]- \vt(P) |\geq 2\delta
\right) \leq
2 \alpha(\e,\delta) e^{-N\gamma(\e,\delta)},
\edeqn
where $f_N^{\lambda_N}\in\F^*_{N,\lambda_N}$.

\end{thm}


 In the literature of machine learning,
consistency analysis refers to (\ref{eq:lemma-saa-approxi-2-a})
whereas in stochastic programming, it refers to 
(\ref{eq:lemma-saa-approxi-2}).
The consistency analysis is mostly focused on
the case when $Z$ is a compact set, we refer readers to  
Norkin and Keyzer \cite{NoK09} which provides an excellent overview about this. 
Caponnetto and Vito \cite{CaV07} is one of a few exceptions
which studies convergence of the empirical risk
of a regularized least-square 
problem in a reproducing kernel Hilbert space with unbounded
feasible set.
Under some moderate conditions,
they derive 
optimal choice of the regularization parameter and optimal rate of convergence of the empirical risk over a class of priors defined by 
a uniformly bounded kernel. 
Our focus here is slightly different: while we are also aiming 
to derive exponential rate of convergence, we concentrate more on how to overcome the complexities and 
challenges arising from a generic form of the cost function and an unbounded kernel. 
For instance, the exponential rate of convergence in (\ref{eq:lemma-saa-approxi-1}) holds uniformly for all $f\in\F$.
This kind of result may not hold in general, see 
a counter example in \cite{SSS10}. Here we
 manage to establish the uniform convergence 
 by showing equi-continuity of the class of functions in $\F$ under Assumption \ref{A:kernel} and their 
 uniform boundedness over a compact subset of $Z$.
 
\noindent
\textbf{Proof of Theorem \ref{thm-vt-consist-ML}.} 
Observe that inequality (\ref{eq:lemma-saa-approxi-1}) implies
\bgeqn
P^{\otimes N} \left(
|\bbe_{P_N}[c(z,f_N^{\lambda_N}(x))]+\lambda_N \|f_N^{\lambda_N}\|_k^2 - \bbe_{P}[c(z,f_N^{\lambda_N}(x))|\geq \delta
\right ) \leq
\alpha(\e,\delta) e^{-N\gamma(\e,\delta)},\nonumber\\
\label{eq:lemma-saa-approxi-1-a}
\edeqn
and a combination of (\ref{eq:lemma-saa-approxi-1-a}) and (\ref{eq:lemma-saa-approxi-2})
yields (\ref{eq:lemma-saa-approxi-2-a}). Thus it suffices to prove (\ref{eq:lemma-saa-approxi-1}) and (\ref{eq:lemma-saa-approxi-2}).
Since $P \in \mathcal{M}_Z^{\phi}$, then for any $\epsilon>0$, there exist
a constant $r>0$  such that
\bgeq
\int_Z \phi(z)\mathbbm{1}_{(r,\infty)}(\phi(z))P (dz) \leq \epsilon.
\edeq
Moreover, by the large deviation theory, there exist positive numbers 
$C_0$ and $\gamma_0$ such that
\bgeq
P^{\otimes N} \left(\int_Z \phi(z)\mathbbm{1}_{(r,\infty)}(\phi(z))P_N (dz) \geq 2\epsilon \right)\leq 
C_0e^{-\gamma_0 N}.
\edeq
Under the coercive condition on $\phi$ in Assumption \ref{A:cost-1} (a),
there exists a compact set $Z_\epsilon=(X_\epsilon,Y_\epsilon)
\subset Z$ such that
$
\{z\in Z: \phi(z) \leq r\} \subset Z_\epsilon.
$
Thus
\bgeqn
\,\,\,\,\sup_{f \in \F}
\int_{Z\backslash Z_\epsilon}
|c(z,f(x))|P(dz)
\leq
\int_{Z\backslash Z_\epsilon}
\phi(z)P(dz)
\leq
\int_{\{z\in Z: \phi(z)> r\}}
\phi(z)P(dz)
 \leq \epsilon
 \label{eq:ZbackslashP}
\edeqn
and
\bgeqn
 \label{eq:ZbackslashPN}
&&P^{\otimes N} \left(
\sup_{f \in \F}
\int_{Z\backslash Z_\epsilon}
|c(z,f(x))|P_N(dz)
\geq  2\epsilon\right)
\leq 
P^{\otimes N} \left(
\int_{Z\backslash Z_\epsilon}
\phi(z)P_N(dz)
\geq 2\epsilon\right)\\
&\leq&
P^{\otimes N} \left(
\int_{\{z\in Z: \phi(z)> r\}}
\phi(z)P_N(dz)
 \geq 2\epsilon \right)
\leq C_0e^{-\gamma_0 N}.\nonumber
\edeqn
By Assumption \ref{A:kernel},
there exists $\eta>0$ such that for 
any $x, x' \in X_\epsilon$ satisfying $\|x-x'\| <\eta$, we have 
\bgeq
|f(x')-f(x)|
&=& |\langle f, k(\cdot,x')\rangle-\langle f, k(\cdot,x)\rangle|\leq \|f\|_k\|k(\cdot,x')-k(\cdot,x)\|_k\nonumber\\
&\leq& \beta \|k(\cdot,x')-k(\cdot,x)\|_k 
\leq \beta \epsilon,
\edeq
which implies $\F$ is equi-continuous when it is restricted to $X_\epsilon$.

Let $\Delta_\epsilon:=\sup_{x \in X_\epsilon}\|k(\cdot,x)\|_k$.
Then for any $f \in \F$,
\bgeq
\sup_{x \in X_\epsilon}|f(x)|
=\sup_{x \in X_\epsilon}|\langle f, k(\cdot,x)| 
\leq \|f\|_k \sup_{x \in X_\epsilon}\|k(\cdot,x)\|_k \leq \beta \Delta_\epsilon,
\edeq
which implies that $\F$ is uniformly bounded when it is restricted to $X_\epsilon$.
Let $\bar{r}:=\max\{|r(z)|: z\in Z_\epsilon$\} and $\bar{\epsilon}:=(\epsilon/{\bar{r}})^{1/\nu}$.
By Ascoli-Arzela Theorem \cite{Br04},
there exists an $\bar{\epsilon}$-net of  $\mathcal{F}_K:=\{f_1,\ldots,f_K\} \subset  \mathcal{F}$
such that
$
\F = \displaystyle \bigcup_{k=1}^K \F_k^{\bar{\epsilon}},
\label{eq:H_K}
$
where
$
\F_k^{\bar{\epsilon}}  := \{f \in \F: \sup_{x\in X_\epsilon}|f(x)-f_k(x)| \leq \bar{\epsilon}\} 
$
for $k=1,\ldots, K$. 
Therefore,
\bgeq
&&|\vt(P_N,\lambda_N)-\vt(P)|\\ 
&=& \left|\sup_{ f\in \F} \{\bbe_{P_N}[c(z,f(x))]+\lambda_N \|f\|_k^2\}
- \sup_{f \in \F} \bbe_{P}[c(z,f(x))]
\right|\\
&\leq& \left|\sup_{ f\in \F} \bbe_{P_N}[c(z,f(x))\mathbbm{1}_{Z_\epsilon}(z))]- \sup_{f \in \F} \bbe_{P}[c(z,f(x))\mathbbm{1}_{Z_\epsilon}(z))]\right|+
\lambda_N\beta^2\\
&&+\sup_{f \in \F}
\int_{Z\backslash Z_\epsilon}
|c(z,f(x))|P_N(dz)+\sup_{f \in \F}
\int_{Z\backslash Z_\epsilon}
|c(z,f(x))|P(dz)\\
&=&\left|
\sup_{k \in K} \sup_{f \in \F_k^{\bar{\epsilon}}} \bbe_{P_N}[c(z,f(x))\mathbbm{1}_{Z_\epsilon}(z))]
-\sup_{k \in K} \sup_{f \in \F_k^{\bar{\epsilon}}} \bbe_{P}[c(z,f(x))\mathbbm{1}_{Z_\epsilon}(z))] \right|+2\epsilon\\
&&+\sup_{f \in \F}
\int_{Z\backslash Z_\epsilon}
|c(z,f(x))|P_N(dz)\\
&\leq& \sup_{k \in \{1,\ldots,K\}} \sup_{f \in \F_k^{\bar{\epsilon}}}
\left|\bbe_{P_N}[c(z,f(x))\mathbbm{1}_{Z_\epsilon}(z))]
-c(z,f_k(x))\mathbbm{1}_{Z_\epsilon}(z))+c(z,f_k(x))\mathbbm{1}_{Z_\epsilon}(z))]\right.\\
&&\left.-\bbe_{P}[c(z,f(x))\mathbbm{1}_{Z_\epsilon}(z))-c(z,f_k(x))\mathbbm{1}_{Z_\epsilon}(z))+c(z,f_k(x))\mathbbm{1}_{Z_\epsilon}(z))]\right|+2\epsilon\\
&&+\sup_{f \in \F}
\int_{Z\backslash Z_\epsilon}
|c(z,f(x))|P_N(dz)\\
&\leq& \sup_{k \in \{1,\ldots,K\}}  \left|\bbe_{P_N}[c(z,f_k(x))\mathbbm{1}_{Z_\epsilon}(z))]- \bbe_{P}[c(z,f_k(x))\mathbbm{1}_{Z_\epsilon}(z))]\right|+4 \epsilon\\
&&+\sup_{f \in \F}
\int_{Z\backslash Z_\epsilon}
|c(z,f(x))|P_N(dz),
\edeq
where the first inequality holds due to 
$\|f\|_k \leq \beta$, and
the last inequality holds  because under
Assumption \ref{A:cost-Hold}
we have
\bgeq
\bbe_{P}[c(z,f(x))\mathbbm{1}_{Z_\epsilon}(z))-c(z,f_k(x))\mathbbm{1}_{Z_\epsilon}(z))]
\leq 
\bbe_{P}[r(z)\|f-f_k\|^{\nu}\mathbbm{1}_{Z_\epsilon}(z)]
\leq \bar{r} \bar{\epsilon}^{\nu}=\epsilon
\edeq
and
\bgeq
\bbe_{P_N}[c(z,f(x))\mathbbm{1}_{Z_\epsilon}(z))-c(z,f_k(x))\mathbbm{1}_{Z_\epsilon}(z))]
\leq 
\bbe_{P_N}[r(z)\|f-f_k\|^{\nu}\mathbbm{1}_{Z_\epsilon}(z)]
\leq \bar{r} \bar{\epsilon}^{\nu}=\epsilon.
\edeq
It follows from 
by the classical Cram\'er's large deviation theorem \cite{DeZ98} 
that for each $k$ there exist positive constants
$C(\e,\delta,f_k)$ and $ \gamma (\e,\delta,f_k)$ such that
\bgeq
 P^{\otimes N} \left(
 |\bbe_{P_N}[c(z,f_k(x))\mathbbm{1}_{Z_\epsilon}(z))]-\bbe_{P}[c(z,f_k(x))\mathbbm{1}_{Z_\epsilon}(z))]| \geq \delta-6 \epsilon \right)
\leq
C(\e,\delta,f_k) e^{-N\gamma (\e,\delta,f_k)}.
\label{eq:approx-phi-phi-N-a}
\edeq
Hence, we have
\bgeq
&&P^{\otimes N} \left(\sup_{f \in {\cal H}} |\bbe_{P_N}[c(z,f(x))]+\lambda_N\|f\|_k^2 - \bbe_{P_N}[c(z,f(x))]| \geq \delta\right)\\
&\leq & P^{\otimes N} \left(\sup_{k \in \{1,\ldots,K\}}  |\bbe_{P_N}[c(z,f_k(x))\mathbbm{1}_{Z_\epsilon}(z))]-\bbe_{P}[c(z,f_k(x))\mathbbm{1}_{Z_\epsilon}(z))]| \geq \delta -6 \epsilon\right) \\
&&+P^{\otimes N} \left(
\sup_{f \in \F}
\int_{Z\backslash Z_\epsilon}
|c(z,f(x))|P_N(dz)
\geq  2\epsilon\right)
\edeq
\bgeq
&\leq& \sum_{k \in \{1,\ldots,K\}} P^{\otimes N} \left(    |\bbe_{P_N}[c(z,f_k(x))\mathbbm{1}_{Z_\epsilon}(z))]-\bbe_{P}[c(z,f_k(x))\mathbbm{1}_{Z_\epsilon}(z))]| \geq \delta -6 \epsilon \right)+C_0e^{-\gamma_0 N}\\
&\leq& \sum_{k \in \{1,\ldots,K\}} C(\e,\delta,f_k) e^{-N\gamma (\e,\delta,f_k)}+C_0e^{-\gamma_0 N},
\edeq
which implies (\ref{eq:lemma-saa-approxi-1}).
\hfill $\Box$

Next we study uniform convergence of the regularized empirical risk with respect to a class of empirical probability distributions as the sample size increases. In practice, 
we may be able to obtain empirical data but often do not know the true probability distribution generating the data. 
Our next result states that the empirical risk converges to its true counterpart uniformly for all empirical data to be used in the machine learning model.

\begin{thm}[Uniform consistency of $\vt(P_N,\lambda_N)$]
\label{thm-vt-consist-ML-uniform}
Let 
Assumptions \ref{A:kernel}, \ref{A:cost-1} and \ref{A:cost-Hold} hold.
Let
$$
\mathcal{M}_{\kappa}^{\phi^p}:=\left\{P\in \mathscr{P}(Z): \int_{Z} \phi(z)^p P(dz)<\kappa\right\}
$$ 
for some fixed $p >1$ and $\mathcal{M}$ be a compact subset of 
$\mathcal{M}_{\kappa}^{\phi^p}$.
 Then for every $\epsilon >0$ and $\delta>0$, there
exists $N_0$ such that  
 \bgeqn
\label{eq:saa-approxi-uni}
\sup_{P \in \mathcal{M}} P^{\otimes N} \left(
|\vt(P_N,\lambda_N) - \vt(P)| \geq \delta
\right) \leq \epsilon,
\edeqn
when $\lambda_N \leq \delta/{4\beta^2}$ and $N \geq N_0$.
\end{thm}

The uniform convergence (\ref{eq:saa-approxi-uni}) is closely related to learnability in statistical learning theory 
which is defined as the uniform convergence of
$R(f_N(P_N))$  to $\vt(P)$ for all empirical probability distributions drawn from
$\mathscr{P}(Z)$, where $R(\cdot)$ is defined as in (\ref{eq:ML-Rf-min}), see
 \cite[Definition 1]{SSS10}. 
 Here we are looking into the convergence for all $P_N$ whose true counterpart  is drawn ${\cal M}$. This applies to the case
that there is some incomplete information about the nature of $P$.

\noindent
\textbf{Proof of Theorem \ref{thm-vt-consist-ML-uniform}.}
We first show that (\ref{eq:saa-approxi-uni}) holds for each $P \in \mathcal{M}  \subset \M^{\phi^p}_{\kappa}$.
For  fixed $\bar{P}$, by the continuity of $\vt(\cdot)$ at $\bar{P}$ in Theorem \ref{t-cont-vt-wrt-P},
for any $\delta >0$, there exists a positive constant $\eta>0$ such that  
\bgeq
|\vt(Q)-\vt(\bar{P})| < \delta/2, 
\edeq
for each $Q$ satisfying $d_{\phi}(Q,\bar{P}) < \eta$.
It follows by \cite[Corollary 3.5]{KSZ12} that
$(\M^{\phi^p}_{\kappa}, \dd_{\phi})$ has the UGC property for all $p> 1$ and $\kappa> 0$,
that is, for any $\epsilon,\eta>0$, there exists $N_0\in \mathbb{N}$ such that for all  $N \geq N_0$
\bgeq
P^{\otimes N}\left[ \dd_{\phi}(P_N,P)\geq \eta\right] \leq \epsilon, \forall P\in \M^{\phi^p}_{\kappa}.
\edeq
Thus, for any $\epsilon >0$ and $\delta>0$, there exists $N_0$ such that for all $N \geq N_0$
\bgeq
\bar{P}^{\otimes N}\left[ |\vt(\bar{P}_N)-\vt(\bar{P})| \geq \delta/2 \right] \leq P^{\otimes N}\left[ \dd_{\phi}(\bar{P}_N,\bar{P}) \geq \eta\right]
\leq \epsilon.
\edeq
Since
\bgeq
|\vt(\bar{P}_N,\lambda_N) - \vt(\bar{P})|
&=&| \inf_{f \in \F} \{\bbe_{\bar{P}}[c(z,f(x))]+\lambda_N \|f\|_k^2\} - \inf_{f \in \F} \bbe_{{\bar{P}}_N}[c(z,f(x))|\\
&\leq& | \inf_{f \in \F} \bbe_{\bar{P}}[c(z,f(x))] - \inf_{f \in \F} \bbe_{{\bar{P}}_N}[c(z,f(x))| +\sup_{f \in \F} \lambda_N \|f\|_k^2 \\
&=&  |\vt({\bar{P}}_N)-\vt(\bar{P})| +\lambda_N \beta^2,
\edeq
then 
\bgeq
 \bar{P}^{\otimes N}\left[ |\vt({\bar{P}}_N,\lambda_N) -\vt(\bar{P})| \geq \delta \right]  \leq  P^{\otimes N}\left[ |\vt(\bar{P}_N)-\vt(\bar{P})| \geq \delta/2 \right] \leq \epsilon
\edeq
when $\lambda_N  \leq \delta/{4\beta^2}$. 
Therefore, (\ref{eq:saa-approxi-uni}) holds when $P$ is  fixed  at $\bar{P}$.

Now we show (\ref{eq:saa-approxi-uni}) holds for all  $P \in \mathcal{M}$.
Assume for the sake of a contradiction that
 there exist some positive numbers $\epsilon_0$ and $\delta_0$ such that for any
$s\in \mathbb{N}$, there exist 
$s'>s$,  $P_{s'}\in  \mathcal{M}$  
and some $N_{s'}\geq s$ such that
\bgeqn
 P_{s'}^{\otimes N_{s'}}\left[|\vt(P_{N_{s'}},\lambda_{N_{s'}}) -\vt(P_{s'})| \geq \delta_0\right] > \epsilon_0. 
 \label{eq:UGC-i-2-a-S-0-a}
\edeqn
Let $s$ increase. Then we obtain a sequence of $\{P_{s'}\}$ which satisfies (\ref{eq:UGC-i-2-a-S-0-a}).
Since $\mathcal{M}$ is compact under the $\phi$-weak topology,
 then $\{P_{s'}\}$ has a converging subsequence.
Assume without loss of generality that
$P_{s'}  \xrightarrow[]{\phi}  P_* \in \mathcal{M}$.
Since $\vt(\cdot)$ is continuous at $P_*$, then there exists $\eta >0$ such that
$|\vt(Q)-\vt(P_*)| < \delta_0/4$ for $P$ satisfying $\dd_{\phi}(Q,P_*) < \eta$  and then 
$$
|\vt(Q,\lambda')-\vt(P_*)| \leq |\vt(Q)-\vt(P_*)| +\lambda' \beta^2< \delta_0/2
$$ 
for $\lambda' \leq \delta_0/{4\beta^2}$.
By $P_{s'}  \xrightarrow[]{\phi}  P_* $, there exists $s'_0$ such that
$\dd_{\phi}(P_{s'}, P_*) < \eta/2$ for $s' \geq s'_0$, 
and then $|\vt(P_{s'},\lambda_{s'})-\vt(P_*)| < \delta_0/2$ for $\lambda_{s'} \leq \delta_0/{4\beta^2}$.
On the other hand, by the UGC property
\bgeq
 P_s^{\otimes N_s} (\dd_\phi(P_{N_{s'}},P_*) \geq \eta) 
 &\leq&  P_s^{\otimes N_s} (\dd_\phi(P_{N_{s'}},P_{s'}) +\dd_\phi(P_{s'},P_*) \geq \eta) \\
 &=&  P_s^{\otimes N_s} (\dd_\phi(P_{N_{s'}},P_{s'})  \geq \eta-\dd_\phi(P_{s'},P_*))\\
  &\leq&  P_s^{\otimes N_s} (\dd_\phi(P_{N_{s'}},P_{s'})  \geq \eta/2)
 \leq \epsilon_0
\edeq
for sufficiently large $N_{s'}$.
Therefore, 
\bgeq
 P_{s'}^{\otimes N_{s'}}\left[|\vt(P_{N_{s'}},\lambda_{N_{s'}}) -\vt(P_*)| \geq \delta_0/2\right] \leq \epsilon_0,
\edeq
and
\bgeq
 &&P_{s'}^{\otimes N_{s'}}\left[|\vt(P_{N_{s'}},\lambda_{N_{s'}}) -\vt(P_{s'})| \geq \delta_0\right]\\
 &\leq&  P_{s'}^{\otimes N_{s'}}\left[|\vt(P_{N_{s'}},\lambda_{N_{s'}}) -\vt(P_*)| + |\vt(P_{s'},\lambda_{s'})-\vt(P_*)| \geq \delta_0\right]\\
  &\leq&  P_{s'}^{\otimes N_{s'}}\left[|\vt(P_{N_{s'}},\lambda_{N_{s'}}) -\vt(P_*)| \geq \delta_0/2\right]
  \leq \epsilon_0,
\edeq
which leads to a contradiction  with (\ref{eq:UGC-i-2-a-S-0-a}) as desired.
\hfill $\Box$

\section{Concluding remarks}

In this paper, we present some theoretical analysis about 
statistical 
robustness of empirical risk in machine learning.
Our focus is on empirical risk but it might be interesting to extend the discussion to kernel learning estimators.
Moreover, our analysis in 
statistical robustness and uniform consistency
does not capture the effect of 
the optimal choice of the regularization parameter in learning process, but we envisage the effect exists and will be helpful to quantify it.
Finally, it might be interesting to carry out some numerical experiments to examine
the statistical robustness of the empirical risk. We leave all these for future research as they require
much more intensive work.


\bibliographystyle{siamplain}
\bibliography{references}

\end{document}